\documentclass[journal]{IEEEtran}
\ifCLASSINFOpdf
\else
\fi

\usepackage[dvipsnames]{xcolor}
\usepackage{mathtools}
\usepackage{cite}
\usepackage{float}
\usepackage{multicol}
\usepackage{multirow}
\usepackage{hyperref}
\usepackage{graphicx}
\usepackage{amsfonts}
\usepackage{amsmath}
\usepackage{color}
\usepackage{epstopdf}
\usepackage{mathtools, nccmath}
\usepackage{multirow}
\usepackage{bbm}
\usepackage{marginnote}
\usepackage{tikz}
\usepackage{hhline}
\usetikzlibrary{arrows,decorations.pathmorphing,backgrounds,fit,positioning,shapes.symbols,chains}
\usepackage{colortbl}
\usepackage{url}
\usepackage{soul}
\usepackage{pgfplots}
\usepackage{tikz}
\usetikzlibrary{shapes.geometric, arrows}
\usetikzlibrary{shapes}
\usepgfplotslibrary{fillbetween}
\usetikzlibrary{patterns}
\usepackage{rotating}
\newenvironment{ldescription}[1]
  {\begin{list}{}%
   {\renewcommand\makelabel[1]{##1\hfill}%
   \settowidth\labelwidth{\makelabel{#1}}%
   \setlength\leftmargin{\labelwidth}
   \addtolength\leftmargin{\labelsep}}}
  {\end{list}}
\usepackage{verbatim}
\usetikzlibrary{shapes.geometric}
\usepackage{circuitikz}
\usepackage{subfigure}
\usepackage{ctable}
\usepackage{enumitem}
\usepackage{dirtytalk}
\definecolor{color_ConvD}{rgb}{0.8588,0.3333,0.3333}
\definecolor{color_bilevel}{rgb}{0.333,0.8588,0.6039}
\definecolor{color_StochD}{rgb}{0.3216,0.6863,0.8980}
\definecolor{mediumaquamarine}{rgb}{0.4, 0.8, 0.67}
\definecolor{moonstoneblue}{rgb}{0.45, 0.66, 0.76}
\definecolor{pastelpurple}{rgb}{0.94, 0.50, 0.50}
\colorlet{color_blue}{blue!50}
\newcommand{\darkercolor}[3]{
    \colorlet{#3}{#1!#2!black}
}
\darkercolor{color_blue}{70}{color_dark_blue}
\definecolor{WildStrawberry}{rgb}{0,0,0} 
\definecolor{ForestGreen}{rgb}{0,0,0} 
\newcommand{\squeezeup}{\vspace{-5.8mm}}
\definecolor{blue}{rgb}{0,0,0}

\begin{document}
\title{Setting Reserve Requirements to Approximate the Efficiency of the Stochastic Dispatch}
\author{Vladimir~Dvorkin~Jr.,~\IEEEmembership{Student member,~IEEE,}
        Stefanos~Delikaraoglou,~\IEEEmembership{Member,~IEEE,}
        and~Juan~M.~Morales,~\IEEEmembership{Senior member,~IEEE}

\thanks{V. Dvorkin Jr. is with the Technical University of Denmark,  Kgs. Lyngby, Denmark (e-mail: vladvo@elektro.dtu.dk).}
\thanks{S. Delikaraoglou is with the ETH Zurich, Zurich, Switzerland (e-mail: delikaraoglou@eeh.ee.ethz.ch).}
\thanks{J.M.Morales is with the University of Malaga, Malaga, Spain (e-mail: juan.morales@uma.es).}
\thanks{
The work by Vladimir Dvorkin Jr. was supported in part by the Russian Foundation for Basic Research (RFBR) according to the research project No. 16-36-00389. The work by Juan M. Morales was supported in part by the Spanish Ministry of Economy, Industry and Competitiveness through project ENE2017-83775-P; and in part by the European Research Council (ERC) under the EU Horizon 2020 research and innovation programme (grant agreement No. 755705) and the Research Program for Young Talented Researchers of the University of Malaga through project PPIT-UMA-B1-2017/18.}
}
\maketitle
\begin{abstract}
This paper deals with the problem of clearing sequential electricity markets under uncertainty. We consider the European approach, where reserves are traded separately from energy to meet exogenous reserve requirements. Recently proposed stochastic dispatch models that co-optimize these services provide the most efficient solution in terms of expected operating costs by computing reserve needs endogenously. However, these models are incompatible with existing market designs. This paper proposes a new method to compute reserve requirements that bring the outcome of sequential markets closer to the stochastic energy and reserves co-optimization in terms of cost efficiency. Our method is based on a stochastic bilevel program that implicitly improves the inter-temporal coordination of energy and reserve markets, but remains compatible with the European market design. We use \textcolor{blue}{two standard IEEE reliability test cases} to illustrate the benefit of intelligently setting operating  reserves in single and multiple reserve control zones.
\end{abstract}


\IEEEpeerreviewmaketitle

\vspace{-0cm}
\section*{Nomenclature}
\textcolor{WildStrawberry}{The main notation used in this paper is stated below. Additional symbols are
defined in the paper where needed.} All symbols are augmented by index $t$ when referring to different time periods.

\vspace{-0cm}
\subsection{Sets and Indices}
\begin{ldescription}{$xxxxx$}
\item [$\Lambda$] Set of transmission lines.
\item [$\omega \in \Omega$] Set of wind power production scenarios.
\item [$i \in I$] Set of conventional generation units.
\item [$j \in J$] Set of loads.
\item [$k \in K$] Set of wind power units.
\item [$n \in N$] Set of nodes.
\item [$z \in Z$] Set of reserve control zones.
\item [$\{\}_{n}$] Mapping of $\{\}$ into the set of nodes.
\item [$\{\}_{z}$] Mapping of $\{\}$ into the set of reserve control zones.
\end{ldescription}

\vspace{-0cm}
\subsection{Decision variables}
\begin{ldescription}{$xxxxx$}
\item [$\delta_{n}^{\text{DA}}$] Day-ahead voltage angle at node $n$ [rad].
\item [$\delta_{n\omega}^{\text{RT}}$] Real-time voltage angle at node $n$ in scenario $\omega$ [rad].
\item [$D_{z}^{\text{U/D}}$] Up-/Downward reserve requirement in zone $z$ [MW].
\item [$L_{j\omega}^{\text{sh}}$] Shedding of load $j$ in scenario $\omega$ [MW].
\item [$P_{i}^{\text{C}}$] Day-ahead \textcolor{WildStrawberry}{dispatch} of conventional unit $i$ [MW].
\item [$P_{k}^{\text{W}}$] Day-ahead \textcolor{WildStrawberry}{dispatch} of wind power unit $k$ [MW].
\item [$P_{k\omega}^{\text{W,sp}}$] Wind spillage of unit $k$ in scenario $\omega$ [MW].
\item [$R_{i}^{\text{U/D}}$] Up-/Downward reserve provision \textcolor{WildStrawberry}{from} unit $i$  [MW].
\item [$r_{i\omega}^{\text{U/D}}$] Up-/Downward reserve deployment of unit $i$ in scenario $\omega$ [MW].
\end{ldescription}

\vspace{-0cm}
\subsection{Parameters}
\begin{ldescription}{$xxxxx$}
\item [$\pi_{\omega}$] Probability of occurrence of wind power production scenario $\omega$.
\item [$C_{i}$] Day-ahead price offer of unit $i$ [\$/MWh].
\item [$C_{i}^{\text{U/D}}$] Up-/Downward reserve price offer of unit $i$ [\$/MW].
\item [$C^{\text{VoLL}}$] Value of lost load [\$/MWh].
\item [$\overline{F}_{nm}$] Capacity of transmission line $(n,m)$ [MW].
\item [$L_{j}$] Demand of load $j$ [MWh].
\item [$\overline{P}_{i}$] Day-ahead quantity offer of unit $i$ [MW].
\item [$\overline{R}_{i}^{\text{U/D}}$] Up-/Downward reserve capacity offer of unit $i$ [MW].
\item [$\widehat{W}_{k}$] Expected generation of wind power unit $k$ [MW].
\item [$W_{k\omega}$] Wind power realization of unit $k$ in scenario $\omega$ [MW].
\item [$X_{nm}$] Reactance of transmission line $(n,m)$ [p.u.].
\end{ldescription}

\vspace{-0cm}
\section{Introduction}

\IEEEPARstart{E}{lectricity} markets are commonly \textcolor{WildStrawberry}{organized in} a sequence of trading floors \textcolor{WildStrawberry}{in which different services are traded in various time-frames.}
\textcolor{WildStrawberry}{According to the European market architecture,}
this sequence consists of reserve and day-ahead markets \textcolor{WildStrawberry}{that are cleared}
12-36 hours \textcolor{WildStrawberry}{before actual power system operation and pertain to trading reserve capacity and energy services, respectively. Getting close to actual delivery of electricity, a real-time market is organized to balance deviations from the initial schedule.} This market design has been established following a conventional view of power system operation, where uncertainty was \textcolor{WildStrawberry}{induced by equipment} contingencies or minor forecast errors of electricity demand. \textcolor{WildStrawberry}{However, considering} the increasing shares of renewable generation, this design has limited ability to cope with variable and uncertain energy sources, while \textcolor{WildStrawberry}{maintaining} a sufficient level of reliability at a reasonable cost \cite{Aigner_2012}.

\textcolor{blue}{To account for the uncertain nature of renewable generation, recent literature proposes economic dispatch models \cite{Morales_2012,1525111}  and unit commitment formulations \cite{Papavasiliou_2015,4806110,4556639} based on stochastic optimization.} Unlike the conventional market design, which downplays the cost of uncertainty, the stochastic model \textcolor{WildStrawberry}{makes use of a probabilistic description of uncertainty and dispatches the system accounting for plausible forecast errors. In this case, reserve requirements are computed endogenously, instead of relying on rule-of-thumb \textcolor{ForestGreen}{methods} such as as the N-1 security criterion. Although the resulting \textit{stochastic ideal} schedule provides} the most efficient solution in terms of expected operating system costs, this design is not adopted in practice due to still unresolved issues like the violation of the least-cost merit-order principle \cite{zavala2017stochastic}.

There are several research contributions devoted to \textit{approximating} the stochastic ideal solution, i.e., approaching the expected operating cost provided by the stochastic dispatch model while sidestepping its theoretical drawbacks, \textcolor{blue}{namely, the violation of \textit{cost recovery} and \textit{revenue adequacy} for certain realizations of the random variables. The cost recovery property guarantees that the  profit of each conventional producer is greater than or equal to its operating costs. The revenue adequacy property requires that the payments that the system operator
must make to and receive from the participants do not cause it to incur a financial
deficit.} \textcolor{WildStrawberry}{Authors} in \cite{morales2014electricity} \textcolor{WildStrawberry}{propose} a new market-clearing procedure \textcolor{WildStrawberry}{according to which wind power is dispatched to a value different than its forecast mean, such that the expected system cost is minimized.  This procedure respects the merit order of the day-ahead market and thus ensures cost recovery of the flexible units.}
An enhanced stochastic dispatch \textcolor{WildStrawberry}{that guarantees} both cost recovery and revenue adequacy \textcolor{WildStrawberry}{for every uncertainty realization} is introduced in \cite{Kazempour_2018}. The main obstacle preventing the implementation of these two models is that they require changing the state of affairs of conventional market structures. \textcolor{blue}{Finally, authors in \cite{SARFATI2018851} propose a stochastic dispatch model that aims at generating proper price signals that incentivize generators to provide reliability services akin to reserves. This model also guarantees cost recovery and revenue adequacy for every uncertainty realization, but in the meantime it does also require significant changes in market design as well as in {\color{blue} the} offering strategies of the renewable power producers.}

\textcolor{WildStrawberry}{More in line with the current practices of the European market design, \cite{Jensen_2017} proposes a systematic method} to adjust available transfer capacities in order to bring operational efficiency of interconnected power systems closer to the stochastic solution.
\textcolor{WildStrawberry}{In the US electricity markets, several Independent System Operators (ISOs), e.g., the California ISO (CAISO) and Midcontinent ISO (MISO) are implementing new ramping capacity products to increase the ramping ability of the system during the real-time re-dispatch in order to cope with steep ramps of net load \cite{Wang_2013}. Essentially, these flexibility products aim to resemble the stochastic dispatch, which inherently finds the optimal allocation of flexible resources between energy and ramping services. } \textcolor{WildStrawberry}{In the same vein, several US ISOs, as for instance the New York ISO, the ISO New England, the MISO, and the Pennsylvania-New Jersey-Maryland (PJM) market, have introduced an operating reserve demand curve (ORDC) in their real-time market \cite{hogan2013electricity}. Motivated by the two-stage stochastic dispatch model, the ORDC mechanism adjusts electricity prices to reflect the scarcity value of reserves for the system operator and incentivize market players to dispatch their units according to a socially optimal schedule.}
\textcolor{WildStrawberry}{The price adjustment through ORDC leads theoretically to perfect arbitrage between energy and reserves in case these two products are co-optimized \cite{papavasiliou2017remuneration}.  However, in the European market that separates energy and reserve capacity trading this arbitrage is inefficient per se, since market players have to value reserves prior to the energy-only market clearing.}

This paper proposes an alternative approach to approximate the stochastic ideal dispatch solution through an intelligent setting of zonal reserve requirements in sequentially cleared electricity markets akin to the European architecture. \textcolor{WildStrawberry}{Here, we solely focus on operating reserves, i.e., generation that is dispatched to respond to net load variations based on economic bids, rather than on regulating services that are activated by automatic generation control.} Traditionally, requirements \textcolor{WildStrawberry}{for operating reserves} are defined based on deterministic security criteria, such as N-1 security constraint violations, where reserves are dimensioned to cover the largest contingency in the system \cite{rebours2005survey}, or based on a mean forecast load error and forced outage rate of system components over a certain horizon, as in the PJM market \cite{manual2012energy}. The main drawback of those approaches is that they ignore the probabilistic nature of renewable generation and neglect the \textcolor{WildStrawberry}{economic} impact of reserve needs on subsequent operations. In order to account for the operational uncertainty, recent literature proposes reserve dimensioning methods based on probabilistic criteria, according to which reserve requirements are drawn from the probabilistic description of uncertainties \textcolor{blue}{ \cite{strbac2007impact,lee2012analyzing,Doherty_1425549,5929570,6942382,7084167,7552596,lange2005uncertainty,5565529,6299425}.}
\textcolor{blue}{
For example, \cite{strbac2007impact} suggests to define the reserve needs such that they cover 97.7\% ($3\sigma$) of the total variation of a Gaussian distribution modeling the joint wind-load uncertainty, disregarding the fact that wind power forecast errors are described by non-Gaussian distributions \cite{lange2005uncertainty}. As a remedy to this drawback, \cite{5565529} proposed a method for setting the reserve requirements using non-parametric probabilistic wind power forecasts. Flying brick and probability box methods in \cite{5929570} and \cite{6942382}, respectively, compute robust envelopes that enclose the net load with a specified probability level. The recent extension of these methods called flexibility envelopes was suggested in \cite{7084167}. These envelopes are based on the same principles but evolve in time to respect the temporal evolution of reserve requirements. As demonstrated in \cite{5929570}, \cite{6942382} and \cite{7552596}, the probabilistic reserve concepts might be integrated into the actual energy management system and derive requirements for capacity, ramping capability and ramping duration of flexible units. In contrast to the deterministic practices,}
the benefit of these methods is that reserve requirements, drawn from accurately predicted distributions, minimize extreme balancing actions provoked by under- or over-procurement of reserves. However, probabilistic requirements are still an exogenous input to the power dispatch, which disregards their potential impact on expected cost.

To this end, we propose a \textcolor{blue}{model} to determine reserves \textcolor{WildStrawberry}{based on a stochastic bilevel programming problem}, which provides the cost-optimal reserve quantities for a European-type market structure. In line with the stochastic \textcolor{WildStrawberry}{dispatch mechanism}, our model computes the reserve requirements that minimize the expected system cost, anticipating their projected \textcolor{WildStrawberry}{impact} on the subsequent operations. Additionally, these requirements are defined accounting for the \textcolor{WildStrawberry}{actual} decision-making process, i.e., the sequence of market-clearing procedures, zonal representation of the power network and the least-cost merit-order principle in all trading floors. As a result, the implementation of these requirements in a conventional market setting, results in a compromise solution between traditional reserve dimensioning practices and the stochastic dispatch model in terms of expected operating cost. Naturally, our approach has limitations: we consider a simplified market setup with a strictly convex representation. Nevertheless, our results do indicate that the intelligent setting of reserve requirements can enhance the short-run cost efficiency of the conventional market with large shares of renewable generation.

\textcolor{blue}{The proposed model can be used as an analytic tool to provide technical and economic insights about the efficacy of different reserve capacity quantification methods, while it can be also used as a decision-support tool by system operators during the reserve setting process. In the latter case, this model can be presumably executed before the day-ahead reserve capacity auction in order to define the reserve requirements that will be used as input in the actual market-clearing process.
Nevertheless, the incorporation of this method in the operational strategy of the system operator does not entail any changes in the existing market setup, since the model output is solely under the discretion of the system operator and  decoupled from market operations.}

The reminder of this paper is organized as follows. Section \ref{sec1} describes the conventional market design and its counterfactual stochastic representation. Section \ref{sec2} introduces the proposed \textcolor{WildStrawberry}{stochastic bilevel programming problem} to compute the optimal reserve requirements that approximate the ideal stochastic solution maintaining the sequential market structure. \textcolor{blue}{Section \ref{section_sol_str} explains the solution strategy based on the multi-cut Bender's algorithm for large-scale applications}. Section \ref{sec3} provides applications of the proposed model to \textcolor{blue}{the IEEE-24 and IEEE-96 reliability test systems}. Section \ref{sec4} concludes the paper.

\section{Electricity market clearing models} \label{sec1}

In this section, we first describe the conventional market structure and the stochastic dispatch model. We then introduce \textcolor{WildStrawberry}{the necessary} modeling assumptions and provide the mathematical formulations of both models.
\subsection{Conventional market and stochastic dispatch framework}

In Europe, power markets are cleared in \textcolor{WildStrawberry}{sequential and independent auctions} which can be represented by the simplified decision-making process illustrated in Fig. \ref{fig:disaptch_det}, which is referred to as the \textit{conventional} market-clearing model.
First, the system operator defines zonal reserve requirements $\mathcal{D}$ based on certain security standards. Then, the reserve capacity market is cleared based on the offer prices and quantities submitted by the flexible producers to find the optimal upward and downward reserve allocation $\Phi^{\text{R}^{*}}$ that minimizes reserve procurement costs $\mathcal{C}^{\text{R}}$.
This allocation accounts for upward and downward reserve requirement constraints included in the set $\mathcal{Q}^{\text{R}}$. At the next stage, power producers submit their price-quantity offers to the day-ahead market that provides the optimal energy schedule $\Phi^{\text{D}^{*}}$ that minimizes the day-ahead energy cost $\mathcal{C}^{\text{D}}$. The set of day-ahead market constraints $\mathcal{Q}^{\text{D}}$ takes into account the reserve capacity $\Phi^{\text{R}^{*}}$ procured at the previous stage. Closer to delivery time, when realization of uncertainty $\omega^\prime$ is known, the system operator runs the real-time market to define a set of optimal re-dispatch actions $\Phi^{\text{B}}_{\omega^\prime}$ that minimizes the balancing cost $\mathcal{C}^{\text{B}}$, considering the previously procured reserve $\Phi^{\text{R}^{*}}$. In this conventional market design, the choice of reserve requirements  $\mathcal{D}$ has a direct impact on the total expected system cost. In fact, the choice of $\mathcal{D}$ influences reserve procurement decisions $\Phi^{\text{R}}$, which in turn affect day-ahead $\Phi^{\text{D}}$ and real-time $\Phi^{\text{B}}$ energy dispatch decisions.

\tikzstyle{Clearing} = [rectangle, rounded corners = 5, minimum width=10, minimum height=10,text centered, draw=black, fill=white!30,line width=0.3mm]
\tikzstyle{Reserve} = [rectangle, rounded corners = 5, minimum width=10, minimum height=10,text centered, draw=black, fill=white!30,line width=0.3mm]

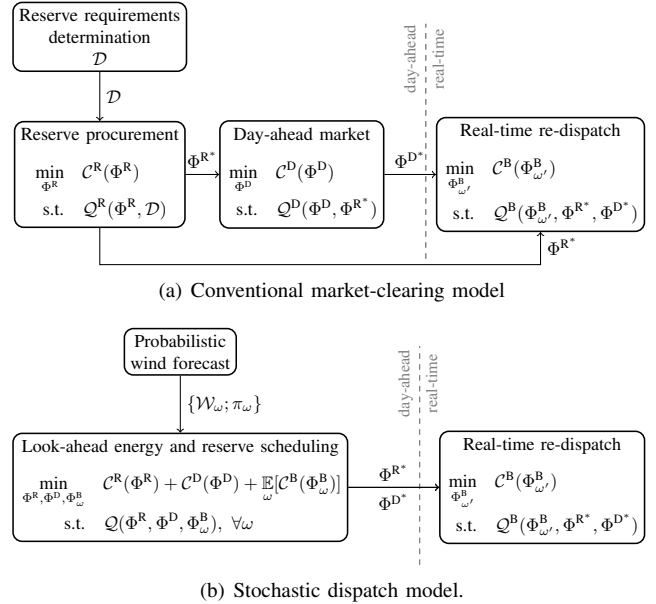
\begin{figure}
\centering
\subfigure[Conventional market-clearing model]
{
\resizebox{8.5cm}{!}{%
\begin{tikzpicture}[node distance=50]
\draw[transform canvas={xshift=177.5, yshift = 0},line width=0.2mm, dashed, color =gray] (0,-4.25) -- node[above, sloped, xshift=45] {\small day-ahead} (0,0.25);
\draw[transform canvas={xshift=177.5, yshift = 0},line width=0.2mm, dashed, color =gray] (0,-4.25) -- node[below, sloped, xshift=45] {\small real-time} (0,0.25);
\node [align=center] (Requirements_box) [Reserve] {Reserve requirements \\ determination \\ $\mathcal{D}$};
\node [align=center] (OR_box) [Clearing, below of = Requirements_box, yshift = -25] {
Reserve procurement  \\ \\
$\begin{aligned}
     \underset{\Phi^{\text{R}}}{\text{min}} \quad& \mathcal{C}^{\text{R}}(\Phi^{\text{R}})\\
     \text{s.t.} \quad& \mathcal{Q}^{\text{R}}(\Phi^{\text{R}}, \mathcal{D})
  \end{aligned}$
};
\draw[transform canvas={xshift=0cm},->,line width=0.2mm] (Requirements_box) -- node[midway, right] {$\mathcal{D}$}  (OR_box) ;
\node [align=center] (DA_box) [Clearing, right of=OR_box, xshift=60] {
Day-ahead market  \\ \\
$\begin{aligned}
     \underset{\Phi^{\text{D}}}{\text{min}} \quad& \mathcal{C}^{\text{D}}(\Phi^{\text{D}})\\
     \text{s.t.} \quad& \mathcal{Q}^{\text{D}} (\Phi^{\text{D}},\Phi^{\text{R}^{*}}) \\
  \end{aligned}$
};
\draw[transform canvas={xshift=0cm},->,line width=0.2mm] (OR_box) -- node[midway, above, sloped] {$\Phi^{\text{R}^{*}}$}  (DA_box) ;
\node [align=center] (RT_box) [Clearing, right of=DA_box, xshift=80] {
Real-time re-dispatch  \\ \\
$\begin{aligned}
     \underset{\Phi_{\omega'}^{\text{B}}}{\text{min}} \quad& \mathcal{C}^{\text{B}}(\Phi_{\omega'}^{\text{B}})\\
     \text{s.t.} \quad& \mathcal{Q}^{\text{B}} (\Phi_{\omega'}^{\text{B}}, \Phi^{\text{R}^{*}},\Phi^{\text{D}^{*}})
  \end{aligned}$
};
\draw[transform canvas={xshift=0cm},->,line width=0.2mm] (DA_box) -- node[midway, above, sloped] {$\Phi^{\text{D}^{*}}$}  (RT_box) ;
\path[draw,->, line width=0.2mm] (OR_box.south) --++ (0,-0.7)  -|   (RT_box.south)  node[midway, right, yshift=10] {$\Phi^{\text{R}^{*}}$} ;
\end{tikzpicture}
}
\label{fig:disaptch_det}
}

\subfigure[Stochastic dispatch model.]
{

\resizebox{8.5cm}{!}{%
\begin{tikzpicture}[node distance=50]
\draw[transform canvas={xshift=140},line width=0.2mm, dashed, color =gray] (5,0.25) -- node[above, sloped, xshift=40] {\small day-ahead} (5,4.5);
\draw[transform canvas={xshift=140},line width=0.2mm, dashed, color =gray] (5,0.25) -- node[below, sloped, xshift=40] {\small real-time} (5,4.5);
\node [align=center] (Forecast_box) [Reserve, xshift=150, yshift=125] {Probabilistic \\ wind forecast};
\node [align=center] (StochD_box) [Clearing, below of = Forecast_box, xshift=0, yshift = -25] {
Look-ahead  energy and reserve \textcolor{WildStrawberry}{scheduling}  \\ \\
$\begin{aligned}
     \underset{\Phi^{\text{R}},\Phi^{\text{D}},\Phi^{\text{B}}_{\omega}}{\text{min}} \quad& \mathcal{C}^{\text{R}}(\Phi^{\text{R}}) + \mathcal{C}^{\text{D}}(\Phi^{\text{D}}) + \underset{\omega}{\mathbb{E}} [  \mathcal{C}^{\text{B}}(\Phi_{\omega}^{\text{B}}) ] \\
     \text{s.t.} \quad& \mathcal{Q}(\Phi^{\text{R}},\Phi^{\text{D}},\Phi^{\text{B}}_{\omega}), \ \forall \omega
  \end{aligned}$
};
\draw[transform canvas={xshift=0cm},->,line width=0.2mm] (Forecast_box) -- node[midway, right] {$\{\mathcal{W}_{\omega}; \pi_{\omega}\}$}  (StochD_box) ;
\node [align=center] (RT_box) [Clearing, right of=StochD_box, xshift=150] {
Real-time re-dispatch  \\ \\
$\begin{aligned}
     \underset{\Phi_{\omega'}^{\text{B}}}{\text{min}} \quad& \mathcal{C}^{\text{B}}(\Phi_{\omega'}^{\text{B}})\\
     \text{s.t.} \quad& \mathcal{Q}^{\text{B}} (\Phi_{\omega'}^{\text{B}}, \Phi^{\text{R}^{*}},\Phi^{\text{D}^{*}})
  \end{aligned}$
};
\draw[transform canvas={xshift=0cm},->,line width=0.2mm] (StochD_box) -- node[midway, above] {$\Phi^{\text{R}^{*}}$}  (RT_box) ;
\draw[transform canvas={xshift=0cm},->,line width=0.2mm] (StochD_box) -- node[midway, below] {$\Phi^{\text{D}^{*}}$}  (RT_box) ;
\end{tikzpicture}
}
\label{fig:disaptch_stoch}
}
\caption{Decision sequences in conventional (a) and stochastic (b) dispatch models.}
\label{fig:sample_subfigures}
\squeezeup
\end{figure}

An alternative model for \textcolor{WildStrawberry}{reserves and energy scheduling} is the \textit{stochastic} \textcolor{WildStrawberry}{dispatch model} outlined in Fig. \ref{fig:disaptch_stoch}. \textcolor{blue}{This is a two-stage stochastic programming model in which first-stage decisions pertain to reserve procurement and day-ahead energy schedule, whereas the second stage models the recourse actions that restore power balance during real-time operation.}
\textcolor{blue}{The stochastic dispatch model} takes as input a probabilistic wind power forecast in the form of a scenario set $\Omega$ and endogenously computes reserve needs. This way, it naturally coordinates all trading floors by co-optimizing reserve ($\Phi^{\text{R}}$) and energy ($\Phi^{\text{D}}$) schedules, anticipating their impact on the subsequent \textcolor{WildStrawberry}{expected balancing cost $\underset{\omega}{\mathbb{E}}  [  \mathcal{C}^{\text{B}}(\Phi_{\omega}^{\text{B}}) ]$ estimated over the scenario set $\Omega$}. \textcolor{blue}{It should be noted that the co-optimization of reserve procurement and energy schedules is a requirement for the implementation of this ideal coordination between the different trading floors.}

In the stochastic dispatch, reserve requirements are a byproduct of the energy and reserve co-optimization problem, resulting in the most efficient solution in terms of total expected operating cost. Moreover, unlike the conventional market model that schedules reserve and day-ahead energy quantities according to the least-cost merit-order principle, the stochastic model schedules \textcolor{WildStrawberry}{generation capacity accounting for potential network congestion during real-time operations, which may lead to} expensive balancing actions \cite{Morales_2012}. This way \textcolor{WildStrawberry}{generators may be scheduled out-of-merit, i.e., more expensive units are dispatched over less expensive ones, in order to minimize the expected costs.}

Despite \textcolor{WildStrawberry}{its superiority in terms of cost efficiency}, the stochastic model suffers from several drawbacks preventing its practical implementation. As already mentioned, the violation of the merit-order principle results in cost recovery and revenue adequacy only in expectation, while for some uncertainty realizations these two essential economic properties may not hold \cite{Morales_2012}.
\textcolor{blue}{This issue disputes the well-functioning of electricity markets in long term, since flexible producers may end up in loss-making positions in one or more scenarios, despite the fact that their expected profit is non-negative. Therefore, these market participants may opt out of the short-run electricity markets or even be discouraged to perform new investments if they are exposed to significant financial risks. In the meantime, the fact that revenue adequacy is only guaranteed in expectation exposes the market operator to the risk of financial deficit. Therefore, a realistic implementation of this market model would require the establishment of out-of-the-market mechanisms, akin to the uplift payments used in the US markets, to provide an ex-post compensation of potential economic deficits. In view of this practical caveats, we do not foresee an actual market clearing implementation of the stochastic dispatch model.}
Moreover, the co-optimization of day-ahead energy and \textcolor{WildStrawberry}{capacity} reserve markets is not compatible with the European market structure, \textcolor{blue}{which dictates that the trading of reserves and energy products is organized in independent sequential auctions.} However, in this work, we \textcolor{WildStrawberry}{show} that the stochastic dispatch solution \textcolor{WildStrawberry}{can} be approximated in the conventional market-clearing model by intelligently setting the reserve requirements $\mathcal{D}$, sidestepping the drawbacks of the stochastic model and improving the efficiency of the \textcolor{WildStrawberry}{existing market setup}.

\subsection{Modeling assumptions}

We use the following set of assumptions to derive computationally tractable yet sensible formulations of the different dispatch models.
Following the European practice, we consider a zonal representation of the network for reserve procurement. In an attempt to build a more generic model, the network topology is included in the day-ahead and real-time dispatch models considering a DC approximation of power flows.
Reserve and energy supply functions are linear, \textcolor{blue}{and all generators are considered to behave as price takers.}
System loads are inelastic with a large value of lost load. This way, the maximization of the social welfare is equivalent to cost minimization.
Flexible units deploy operating reserves with marginal costs of production. The incentive to provide flexibility services is accounted for in reserve offering prices.
\textcolor{blue}{Following the prevailing portfolio bidding adopted in the European markets \cite{6487424}, we consider that all unit commitment and inter-temporal constraints are integrated into the bidding strategies of the generating units. For instance, the commitment of thermal units in practice might be controlled by market participants when offering at either zero price or market price cap. Similarly, offering a part of capacity at zero and even negative price ensures the compliance with the technical minimum constraint of thermal units.}
\textcolor{blue}{This approach is compatible with the European market structure and preserves the convexity of the reserve capacity and day-ahead market-clearing algorithms. In principle, the proposed model can be also applied to market designs that involve non-convex constraints, as for instance the majority of electricity markets in the US, using tight convex relaxations of the unit commitment binary variables. However, this approach lies out of scope of this paper, but we refer the interested reader to \cite{Kasina_2014,7914790} for further discussion}.
Finally, uncertainty is described by a finite set of scenarios and solely induced by stochastic wind power production.

\subsection{Mathematical formulation}

\subsubsection{Conventional market-clearing model}
The sequential procedure, sketched in \textcolor{WildStrawberry}{Fig.} \ref{fig:disaptch_det}, for each hour of the next day is modeled by the following three linear optimization problems.

The reserve procurement problem writes as:
\begingroup
\allowdisplaybreaks
\begin{subequations} \label{prob:reserve_clearing}
\begin{align}
\underset{\Xi^{\text{OR}}}{\text{min}} \quad&
\sum_{i \in I}
\Big(C_{i}^{\text{U}}  R_{i}^{\text{U}} + C_{i}^{\text{D}}  R_{i}^{\text{D}}\Big)  \label{objRC}\\
\text{s.t.} \quad&\sum_{i \in I_{\textcolor{WildStrawberry}{z}}}R_{i}^{\text{U}} = D_{z}^{\text{U}}, \quad \sum_{i \in I_{\textcolor{WildStrawberry}{z}} }R_{i}^{\text{D}} = D_{z}^{\text{D}}, \quad \forall z \in Z, \label{RC:demand}\\
& R_{i}^{\text{U}} + R_{i}^{\text{D}} \leq \overline{P}_{i}, \quad \forall i \in I, \label{RC:limits_UD} \\
&0 \leq R_{i}^{\text{U}} \leq \overline{R}_{i}^{\text{U}}, \quad 0 \leq R_{i}^{\text{D}} \leq \overline{R}_{i}^{\text{D}}, \quad \forall i \in I, \label{RC:limits}
\end{align}
\end{subequations}
\endgroup
where $\Xi^{\text{OR}} = \{R_{i}^{\text{U}}, R_{i}^{\text{D}}, \textcolor{WildStrawberry}{\forall i}\}$ is the set of optimization variables comprising the upward and downward reserve \textcolor{WildStrawberry}{schedule} per each flexible generator. Optimal $\Xi^{\text{OR*}}$ minimizes the reserve procurement cost given by (\ref{objRC}). Equality constraints (\ref{RC:demand}) ensure that zonal reserve upward and downward requirements, denoted as $D_{z}^{\text{U}}$ and $D_{z}^{\text{D}}$, respectively, are fulfilled, whereas inequality constraints (\ref{RC:limits_UD}) - (\ref{RC:limits}) account for the quantity offers of each flexible generator.

Once reserve allocation $\{R_{i}^{\text{U}*}, R_{i}^{\text{D}*}, \forall i\}$ is determined, the least-cost day-ahead energy schedule is computed solving the following optimization problem:
\begingroup
\allowdisplaybreaks
\begin{subequations} \label{prob:day_ahead_clearing}
\begin{align}
\underset{\Xi^{\text{DA}}}{\text{min}} \quad& \sum_{i \in I} C_{i}  P_{i}^{\text{C}}  \label{obj:DA}\\
\text{s.t.} \quad &\sum_{i \in I_{n}}P_{i}^{\text{C}} + \sum_{k \in K_{n}}P_{k}^{\text{W}} - \sum_{j \in J_{n}}L_{j} \nonumber\\
& - \sum_{m:(n,m)\in\Lambda}\frac{\delta_{n}^{\text{DA}}-\delta_{m}^{\text{DA}}}{x_{nm}} = 0,  \quad \forall n \in N,  \label{DA:balance} \\
&R_{i}^{\text{D}*} \leq P_{i}^{\text{C}} \leq \overline{P}_{i} - R_{i}^{\text{U*}}, \quad \forall i \in I, \label{DA:conv_cap}\\
&0 \leq P_{k}^{\text{W}} \leq \widehat{W}_{k}, \quad \forall k\in K, \label{DA:wind_cap}\\
&\frac{\delta_{n}^{\text{DA}}-\delta_{m}^{\text{DA}}}{x_{nm}} \leq \overline{F}_{nm}, \quad \forall (n,m) \in \Lambda,  \label{DA:flow_cap}
\end{align}
\end{subequations}
\endgroup
where $\Xi^{\text{DA}} = \{P_{i}^{\text{C}}, \textcolor{WildStrawberry}{\forall i}; P_{k}^{\text{W}}, \textcolor{WildStrawberry}{\forall k}; \delta_{n}^{\text{DA}}, \textcolor{WildStrawberry}{\forall n}\}$ is the set of variables including day-ahead energy quantities for each conventional and stochastic generator as well as voltage angles at each node. The objective function (\ref{obj:DA}) to be minimized is the day-ahead energy cost, subject to nodal power balance constraints (\ref{DA:balance}), offering limits of conventional and stochastic generators (\ref{DA:conv_cap})-(\ref{DA:wind_cap}) and transmission capacity limits (\ref{DA:flow_cap}). Note that the reserve procurement decisions from the previous stage limit the dispatch of flexible generators at the day-ahead stage. In this design, stochastic production is bounded by the conditional expectation $\widehat{W}_{k}$.

\textcolor{WildStrawberry}{Getting closer to real-time operation, any deviation from the optimal day-ahead dispatch $\{P_{i}^{\text{C}*}, \forall i; P_{k}^{\text{W}*}, \forall k; \delta_{n}^{\text{DA}*}, \forall n\}$ has to be covered by proper balancing actions.} For a specific realization of stochastic production \textcolor{WildStrawberry}{$W_{k\omega'}$}, the optimal re-dispatch is found solving the following linear programming problem:
\begingroup
\allowdisplaybreaks
\begin{subequations} \label{prob:real_time_clearing}
\begin{align}
\underset{\Xi^{\text{RT}}}{\text{min}} & \quad\sum_{i \in I}  C_{i} \Big(r_{i\omega'}^{\text{U}}-r_{i\omega'}^{\text{D}}\Big) +\sum_{j \in J} C^{\text{VoLL}} L_{j\omega'}^{\text{sh}} \label{objRT} \\
\text{s.t.}& \quad \sum_{i \in I_{n}} \left( r_{i\omega'}^{\text{U}}-r_{i\omega'}^{\text{D}} \right)
+ \sum_{k \in K_{n}}\left( W_{k\omega'} - P_{k}^{\text{W*}} - P_{k\omega'}^{\text{W,sp}} \right) \nonumber \\
&  +\sum_{j \in J_{n}}L_{j\omega'}^{\text{sh}} - \!\!\!
  \sum_{m:(n,m)\in\Lambda} \!\!\!\! \frac{\delta_{n\omega'}^{\text{RT}}-\delta_{n}^{\text{DA*}}-\delta_{m\omega'}^{\text{RT}}+\delta_{m}^{\text{DA*}}}{x_{nm}} \nonumber\\
& = 0, \quad \forall n \in N, \label{RT:balance}\\
&0 \leq r_{i\omega'}^{\text{U}} \leq R_{i}^{\text{U}*}, \quad 0 \leq r_{i\omega'}^{\text{D}} \leq R_{i}^{\text{D}*}, \quad \forall i \in I, \label{RT: updownlim}\\
&\frac{\delta_{n\omega'}^{\text{RT}}-\delta_{m\omega'}^{\text{RT}}}{x_{nm}} \leq \overline{F}_{nm}, \quad \forall (n,m) \in \Lambda, \label{RT: maxcapline}\\
&0 \leq P_{k\omega'}^{\text{W,sp}} \leq W_{k\omega'}, \quad \forall k \in K, \label{RT: spill}\\
&0 \leq L_{j\omega'}^{\text{sh}} \leq L_{j}, \quad \forall j \in J, \label{RT: shed}
\end{align}
\end{subequations}
\endgroup
where $\Xi^{\text{RT}} = \{r_{i\omega'}^{\text{U}}, r_{i\omega'}^{\text{D}}, \textcolor{WildStrawberry}{\forall i}; L_{j\omega'}^{\text{sh}}, \textcolor{WildStrawberry}{\forall j}; P_{k\omega'}^{\text{W,sp}}, \textcolor{WildStrawberry}{\forall k}; \delta_{n\omega'}^{\text{RT}}, \textcolor{WildStrawberry}{\forall n}\}$ is the set of re-dispatch decisions, comprising activation of operating reserves, load shedding, wind spillage and real-time voltage angles. The objective function (\ref{objRT}) to be minimized is the balancing cost. Equality constraints (\ref{RT:balance}) ensure the real-time nodal power balance. Inequalities (\ref{RT: updownlim}) limit activation of upward and downward reserves considering the procured reserve quantities. Constraints (\ref{RT: maxcapline}) account for the power capacity of transmission lines. Finally, inequalities (\ref{RT: spill}) and (\ref{RT: shed}) limit wind spillage and load shedding actions to the actual realization of production and system demand, respectively.

\subsubsection{Stochastic dispatch model}
Assuming that wind power uncertainty is described by a finite set of outcomes $W_{k\omega}$ with corresponding probabilities  $\pi_{\omega}$, the stochastic dispatch procedure outlined in Fig. \ref{fig:disaptch_stoch} writes as follows:
\begingroup
\allowdisplaybreaks
\begin{subequations} \label{stochastic_dis}
\begin{align}
\underset{\Xi^{\text{SD}}}{\text{min}} \quad & \sum_{i \in I}
\Big(C_{i}^{\text{U}}  R_{i}^{\text{U}} + C_{i}^{\text{D}}  R_{i}^{\text{D}} + C_{i}  P_{i}^{\text{C}} \Big) + \nonumber\\
& \sum_{\omega} \pi_{\omega} \Big( \sum_{i \in I}  C_{i} \Big(r_{i\omega}^{\text{U}}-r_{i\omega}^{\text{D}}\Big) +\sum_{j \in J} C^{\text{VoLL}} L_{j\omega}^{\text{sh}} \Big) \label{SD_obj}\\
\text{s.t.} \quad& \text{constraints \eqref{RC:demand} - \eqref{RC:limits}} \label{stochD_RCc}\\
& \text{constraints \eqref{DA:balance} - \eqref{DA:flow_cap}} \label{stochD_DAc}\\
& \text{constraints \eqref{RT:balance} - \eqref{RT: shed}}, \quad \forall \omega \in \Omega \label{stochD_RTc}
\end{align}
\end{subequations}
\endgroup
where $\Xi^{\text{SD}} = \{\Xi^{\text{OR}} \cup \Xi^{\text{DA}} \cup \Xi^{\text{RT}}, \textcolor{WildStrawberry}{\forall \omega} \cup (D^{\text{U}},D^{\text{D}})\}$ is the set of stochastic dispatch variables. The objective function (\ref{SD_obj}) to be minimized is the reserve and day-ahead energy cost as well as the expectation of the real-time cost, i.e., the expected cost over the entire decision sequence. Note, that upward and downward reserve \textcolor{WildStrawberry}{requirements  $D_{z}^{\text{U}}$ and $D_{z}^{\text{D}}$ in (\ref{RC:demand}) are decision variables} and only used to reveal \textcolor{WildStrawberry}{optimal} reserve requirements \textcolor{WildStrawberry}{in a stochastic programming sense.}

After the optimal reserve procurement and day-ahead energy schedule are obtained, the system operator solves the real-time re-dispatch problem for a specific realization of the stochastic production $\omega'$ using formulation (\ref{prob:real_time_clearing}).

\section{Approximating the stochastic ideal} \label{sec2}

On the one hand, the conventional procedure has limited capability to accommodate large shares of stochastic production in a cost efficient manner compared to the stochastic dispatch. On the other hand, the adoption of the stochastic procedure appears to be unrealistic because \textcolor{blue}{it does not guarantee revenue adequacy and cost recovery for every uncertainty realization; these are important properties that, in contrast, hold in the sequential market structure \cite{Morales_2012,morales2014electricity}}. For this reason, our motivation is to enhance the cost-efficiency of the conventional market-clearing procedure without changing the market structure. In this line, we introduce a model that approximates the ideal stochastic solution  within the conventional dispatch model by the appropriate setting of zonal reserve requirements. \textcolor{WildStrawberry}{ In essence, we aim at finding the reserve requirements that plugged into the conventional market-clearing model (\ref{prob:reserve_clearing})-(\ref{prob:real_time_clearing}) will yield the minimum total expected system cost.} To compute them, we use the \textcolor{WildStrawberry}{bilevel programming problem} illustrated in Fig. \ref{fig:optimal_req_determination}.

This model comprises two levels. The objective function of the upper level is the same as \eqref{SD_obj} in the stochastic model \eqref{stochastic_dis} and aims at minimizing the total expected system cost. The upper-level constraints enforce real-time re-dispatch limits.
The lower level consists of two optimization problems, \textcolor{WildStrawberry}{namely, the reserve procurement and day-ahead market clearing problems, which are identical to the corresponding optimization problems \eqref{prob:reserve_clearing} and \eqref{prob:day_ahead_clearing} of the conventional model. However, in this bilevel structure, reserve requirements $\mathcal{D}$ are decision variables of the upper-level problem, entering as parameters in the lower-level reserve procurement problem. Hence, reserve requirements $\mathcal{D}$ are not an exogenous input to this model but are internally optimized, accounting for their impact in all three trading floors. As shown in Fig. \ref{fig:optimal_req_determination}, the upper-level decision on $\mathcal{D}$ affects the reserve procurement schedule in the first lower-level problem, which in turn impacts the day-ahead clearing obtained from the second lower-level problem.} In addition, the reserve and energy schedules $\Phi^{\text{R}}$ and $\Phi^{\text{D}}$ enter the upper level, constraining the real-time re-dispatch decisions.

\textcolor{WildStrawberry}{The structure of this stochastic bilevel model guarantees that the temporal sequence of the different markets follows the existing European paradigm.
Having the reserve capacity and day-ahead market clearings as two independent lower-level problems, ensures that reserves and day-ahead schedules are optimized separately, i.e., there is no co-optimization of energy and reserves, while none of these markets have information about the future re-dispatch actions. This property suffices to reproduce the real-time re-dispatch for each scenario independently by including the corresponding constraints only in the upper-level problem. }

Compared to the stochastic model, the main advantage of this bilevel scheme is that it respects the merit-order principle in the reserve capacity and day-ahead energy markets. In fact, given the same reserve requirements, the solutions of both lower-level problems are identical to the solutions of problems \eqref{prob:reserve_clearing} and \eqref{prob:day_ahead_clearing}. \textcolor{WildStrawberry}{Nonetheless, the upper-level problem can still anticipate the impact of reserve requirements on all trading floors and consequently on the total expected cost.}

Since this model is solved prior to any market-clearing procedure, we assume that the system operator can gather information on the price-quantity offers of market participants.
\textcolor{WildStrawberry}{Even in the case of having to use an estimation of price-quantity offers similar to the ORDC mechanism, our approach accounts systematically for the impact of reserve procurement and the structure of forecast errors in all three trading floors.}
In a more realistic setup, this \textcolor{WildStrawberry}{information can be obtained} using inverse optimization techniques as proposed in \cite{ruiz2013revealing} and \cite{mitridati2017bayesian}.

\textcolor{WildStrawberry}{Mathematically, the proposed} reserve determination model writes as the following stochastic bilevel programming problem:
\begingroup
\allowdisplaybreaks
\begin{subequations} \label{prob:bilevel_clearing}
\begin{align}
&\underset{\Xi^{\text{RT}}, D_{z}^{\text{U}}, D_{z}^{\text{D}}}{\text{min}} \quad \eqref{SD_obj} \label{Bilevel_obj}\\
& \;\;\; \text{s.t.} \;\; \text{constraints \eqref{RT:balance} - \eqref{RT: shed}}, \quad \forall \omega \in \Omega,  \label{bilevelcon1}\\
& \;\;\;  D_{z}^{\text{U}}, D_{z}^{\text{D}} \geq 0, \quad \forall z \in Z, \label{bilevelcon2}\\
& \;\;\; ( R_{i}^{\text{U}}, R_{i}^{\text{D}} ) \in \text{arg}
	\left\{\!\begin{aligned}
	&  \underset{\Xi^{\text{OR}}}{\text{min}} \quad \eqref{objRC} \\
	&  \text{s.t.} \;\; \text{constraints \eqref{RC:demand} - \eqref{RC:limits}}
	\end{aligned}\right\}, \label{eq:LLR} \\
	& \;\; \left(\begin{subarray}{c} P_{i}^{\text{C}}, P_{k}^{\text{W}}, \\  \delta_{n}^{\text{DA}} \end{subarray} \right)  \in \text{arg}
	\left\{\!\begin{aligned}
	&  \underset{\Xi^{\text{DA}}}{\text{min}} \quad \eqref{obj:DA} \\
	&  \text{s.t.} \;\; \text{constraints \eqref{DA:balance} - \eqref{DA:flow_cap}}
	\end{aligned}\right\}. \label{eq:LLD}
\end{align}
\end{subequations}
\endgroup

\textcolor{blue}{
According to the mathematical structure of model \eqref{prob:bilevel_clearing}, the lower-level problems \eqref{eq:LLR} and \eqref{eq:LLD} guarantee that the reserve capacity and day-ahead energy markets are serially and independently optimized. This property is in accordance with the time-line of these trading floors in the European market framework. This temporal sequence is accomplished considering that upward $R_i^\text{U*}$ and downward $R_i^\text{D*}$ reserve schedules are variables of the reserve capacity market \eqref{eq:LLR} but enter as parameters in the day-ahead energy market \eqref{eq:LLD}. Moreover, neither problem \eqref{eq:LLR} nor \eqref{eq:LLD} can foresee the outcome of the balancing market, which is included in the upper level of model \eqref{prob:bilevel_clearing}. As a result, both markets have no information about the effect of their decisions on the real-time market. In turn, constraints \eqref{bilevelcon1}-\eqref{bilevelcon2} and the third term of the objective function \eqref{SD_obj} clear the real-time market of the conventional model (\ref{prob:reserve_clearing})-(\ref{prob:real_time_clearing}), independently for each scenario $\omega \in \Omega$, considering that the real-time re-dispatch cannot impact the previous trading floors which are `fixed' to the conventional market solution through the lower-level problems \eqref{eq:LLR} and \eqref{eq:LLD}.}

\tikzstyle{Clearing} = [rectangle, rounded corners = 5, minimum width=10, minimum height=10,text centered, draw=black, fill=white!30,line width=0.3mm]
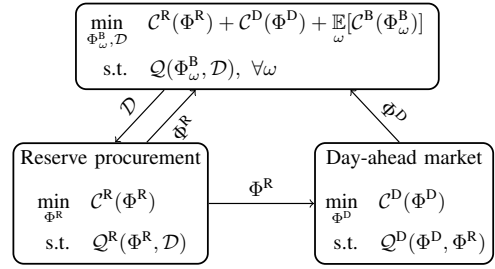
\begin{figure}[]
\center
\resizebox{6.5cm}{!}{%
\begin{tikzpicture}[node distance=50]
\node [align=center] (Upper_level) [Clearing, xshift=0, yshift = -25] {
$\begin{aligned}
     \underset{\Phi^{\text{B}}_{\omega}, \mathcal{D}}{\text{min}} \quad& \mathcal{C}^{\text{R}}(\Phi^{\text{R}}) + \mathcal{C}^{\text{D}}(\Phi^{\text{D}}) + \underset{\omega}{\mathbb{E}} [  \mathcal{C}^{\text{B}}(\Phi_{\omega}^{\text{B}}) ] \\
     \text{s.t.} \quad& \mathcal{Q}(\Phi^{\text{B}}_{\omega}, \mathcal{D}), \ \forall \omega
  \end{aligned}$
};
\node [align=center] (OR_box) [Clearing, below of = Upper_level, yshift = -25, xshift = -70] {
Reserve procurement  \\ \\
$\begin{aligned}
     \underset{\Phi^{\text{R}}}{\text{min}} \quad& \mathcal{C}^{\text{R}}(\Phi^{\text{R}})\\
     \text{s.t.} \quad& \mathcal{Q}^{\text{R}}(\Phi^{\text{R}}, \mathcal{D})
  \end{aligned}$
};
\node [align=center] (DA_box) [Clearing, below of = Upper_level, yshift = -25, xshift = 70] {
Day-ahead market  \\ \\
$\begin{aligned}
     \underset{\Phi^{\text{D}}}{\text{min}} \quad& \mathcal{C}^{\text{D}}(\Phi^{\text{D}})\\
     \text{s.t.} \quad& \mathcal{Q}^{\text{D}} (\Phi^{\text{D}},\Phi^{\text{R}}) \\
  \end{aligned}$
};
\draw[transform canvas={xshift=-25},->,line width=0.2mm] (Upper_level) -- node[midway, above, sloped] {$\mathcal{D}$}  (OR_box) ;
\draw[transform canvas={xshift=-10},->,line width=0.2mm] (OR_box) -- node[midway, below, sloped] {$\Phi^{\text{R}}$}  (Upper_level) ;
\draw[transform canvas={xshift=0},->,line width=0.2mm] (OR_box) -- node[midway, above, sloped] {$\Phi^{\text{R}}$}  (DA_box) ;
\draw[transform canvas={xshift=25},->,line width=0.2mm] (DA_box) -- node[midway, above, sloped, xshift=5] {$\Phi^{\text{D}}$}  (Upper_level) ;
\end{tikzpicture}
}
\caption{Bilevel structure of the proposed reserve determination model.}
\label{fig:optimal_req_determination}
\squeezeup
\end{figure}

This formulation is \textcolor{WildStrawberry}{computationally} intractable, since it consists of an upper-level \textcolor{WildStrawberry}{optimization problem} constrained by two lower-level optimization problems. However, since both lower-level problems are convex with linear objective functions and constraints, they can be replaced by their Karush–-Kuhn–-Tucker optimally conditions, such that the problem can be recast as as a single-level \textcolor{WildStrawberry}{mathematical program with equilibrium constraints} (MPEC). The resulting model \textcolor{WildStrawberry}{includes a set of}  nonlinear \textcolor{WildStrawberry}{complementary} slackness constraints, which can be linearized using \textcolor{WildStrawberry}{disjunctive constraints} or SOS1 variables, transforming the MPEC problem into a \textcolor{WildStrawberry}{mixed-integer linear program} (MILP) \cite{pozo2017basic}.

\textcolor{blue}{
\section{Solution strategy} \label{section_sol_str}
The set of integer variables used to linearize the complementarity constraints of the lower-level problems \eqref{eq:LLR} and \eqref{eq:LLD} limits the application of the proposed reserve quantification model to power systems of moderate scale. For the large-scale applications, we propose an iterative solution strategy based on the multi-cut Bender's algorithm \cite{conejo2006decomposition}. For a fixed reserve and day-ahead dispatch, the set of real-time constraints \eqref{RT:balance} - \eqref{RT: shed} is independent per scenario. This allows for Bender's decomposition where each subproblem solves a scenario-specific real-time re-dispatch problem. The subproblems at iteration $\nu$ write as follows:
\begingroup
\allowdisplaybreaks
\begin{subequations} \label{benders_sub_asd}
\begin{align}
\Bigg\{\underset{\Xi_{s}^{\text{RT,B}}}{\text{min}} \quad & C^{\text{RT}(\nu)}_{\omega} := \sum_{i \in I}  C_{i} \Big(r_{i\omega}^{\text{U}}-r_{i\omega}^{\text{D}}\Big) +\sum_{j \in J} C^{\text{VoLL}} L_{j\omega}^{\text{sh}} \label{objRTbenderssub} \\
\text{s.t.} \quad & R_{i}^{\text{U}} = \tilde{R}_{i}^{\text{U}(\nu)} \quad: \theta_{i\omega}^{R_{i}^{\text{U}}(\nu)}, \quad \forall i \in I, \label{bend_sub_1}\\
& R_{i}^{\text{D}} = \tilde{R}_{i}^{\text{D}(\nu)} \quad: \theta_{i\omega}^{R_{i}^{\text{D}}(\nu)}, \quad \forall i \in I,\\
& P_{k}^{\text{W}} = \tilde{P}_{k}^{\text{W}(\nu)} \quad: \theta_{k\omega}^{P_{k}^{\text{W}}(\nu)}, \quad \forall k \in K,\\
& \delta_{n}^{\text{DA}} = \tilde{\delta}_{n}^{\text{DA}(\nu)} \quad: \theta_{n\omega}^{\delta_{n}^{\text{DA}}(\nu)}, \quad \forall n \in N, \label{bend_sub_4}\\
& \text{constraints \eqref{RT:balance} - \eqref{RT: shed}} \Bigg\} \quad \forall \omega \in \Omega, \nonumber
\end{align}
\end{subequations}
\endgroup
where $\Xi_{s}^{\text{RT,B}} = \Xi^{\text{RT}} \cup \{
R_{i}^{\text{U}}, R_{i}^{\text{D}}, \forall i; P_{k}^{\text{W}}, \forall k; \delta_{n}^{\text{DA}}, \forall n
\}$ is the set of decision variables of each subproblem of the Bender's algorithm. Constraints \eqref{bend_sub_1} - \eqref{bend_sub_4} fix the first-stage decisions to their optimal values obtained at the previous iteration, and the corresponding dual variables yield sensitivities of the reserve and day-ahead decisions used in Bender's cuts.
}

\textcolor{blue}{
The master problem of the Bender's algorithm at iteration $\nu$ writes as follows:
\begingroup
\allowdisplaybreaks
\begin{subequations} \label{prob:real_time_clearing_benders}
\begin{align}
\underset{\Xi^{\text{M,B}}}{\text{min}} \quad&  \sum_{i \in I}
\Big(C_{i}^{\text{U}}  R_{i}^{\text{U}} + C_{i}^{\text{D}}  R_{i}^{\text{D}} + C_{i}  P_{i}^{\text{C}} \Big)
+ \sum_{\omega \in \Omega} \pi_{\omega} \alpha_{\omega}^{(\nu)} \\
\text{s.t.} \quad & \alpha_{\omega}^{(\nu)} \geq C^{\text{RT}(\rho)}_{\omega} + \sum_{i \in I} \theta_{i\omega}^{R_{i}^{\text{U}}(\rho)} \Big(R_{i}^{\text{U}} - R_{i}^{\text{U}(\rho)} \Big) \nonumber\\
& \quad\quad + \sum_{i \in I} \theta_{i\omega}^{R_{i}^{\text{D}}(\rho)} \Big(R_{i}^{\text{D}} - R_{i}^{\text{D}(\rho)} \Big) \nonumber \\
& \quad\quad + \sum_{k \in K} \theta_{k\omega}^{P_{k}^{\text{W}}(\rho)} \Big(P_{k}^{\text{W}} - P_{k}^{\text{W}(\rho)} \Big) \nonumber \\
& \quad\quad + \sum_{n \in N} \theta_{n\omega}^{\delta_{n}^{\text{DA}}(\rho)} \Big(\delta_{n}^{\text{DA}} - \delta_{n}^{\text{DA}(\rho)} \Big), \nonumber \\
&\quad\quad\quad \rho = 1 \dots \nu-1, \forall \omega \in \Omega, \label{ben_cut} \\
&\alpha_{\omega}^{(\nu)} \geq \underline{\alpha}, \quad \forall \omega \in \Omega, \label{ben_cut_zero}\\
&D_{z}^{\text{U}}, D_{z}^{\text{D}} \geq 0, \quad \forall z \in Z, \\
&\text{Linearized KKT conditions of \eqref{eq:LLR}}, \\
&\text{Linearized KKT conditions of \eqref{eq:LLD}},
\end{align}
\end{subequations}
\endgroup
where $\Xi^{\text{M,B}} = \Xi^{\text{OR}} \cup \Xi^{\text{DA}} \cup \alpha_{\omega}$ is the set of decisions variables of the master problem, and index $\rho$ is used to integrate the fixed values of the corresponding variables at previous iterations. The Bender's cuts are updated at each iteration by \eqref{ben_cut} using sensitivities from all previous iterations, while \eqref{ben_cut_zero} imposes a lower bound $\underline{\alpha}$ on the auxiliary variable $\alpha$.  Since the subproblems allow for load shedding, they are always feasible, requiring no feasibility cuts in the master problem. The algorithm converges at iteration $\nu$ if $\Big|\sum_{\omega \in \Omega} \pi_{\omega} \big(\alpha_{\omega}^{(\nu)}-C^{\text{RT}(\nu)}_{\omega}\big)\Big| \leq \epsilon$, where $\epsilon$ is a predefined tolerance.
}

\section{Case Study} \label{sec3}
In this section, we first describe the test system in \textcolor{blue}{Section} \ref{Desription_of_test_system}. In \textcolor{blue}{Section}  \ref{Impact_of_reserve_requirements} and \textcolor{blue}{Section}  \ref{Approximating_the_stoch} we study the impact of reserve requirements on expected operating costs and we assess the remaining efficiency gap of our model with respect to the stochastic solution for a single reserve control zone. In \textcolor{blue}{Section}  \ref{Optimal_zonal_reserve} we extend our analysis to the case of multiple reserve control zones. \textcolor{blue}{In Section \ref{24RTS_with_UC_conctrains} we assess the model's performance in the presence of non-convex technical constraints. Finally, in Section \ref{IEEE96study} we demonstrate the scalability of the model using the proposed Bender's decomposition algorithm.
}

\begin{figure*}
    \centering
    \subfigure[Reserve procurement cost]
    {
\resizebox{8cm}{!}{%
\begin{tikzpicture}[thick,scale=1]
\begin{axis}[view={0}{90},
colormap/viridis,
        xmin=50,
        xmax=350,
        ymax=160,
        ymin=30,
        colorbar,
        colorbar style=     {
        yticklabel style=   {
        try min ticks=6,
        /pgf/number format/.cd,
            fixed,
            fixed zerofill,
            precision=1,
        /tikz/.cd
                            },
        ylabel={Cost [\$1000]},
        y label style=  {
        rotate=0,
        yshift = -80,
                        },
        ymin=0.253,
        ymax=2.0834000,
                            },
        ylabel={$D^{\text{D}}$ [MW]},
        xlabel={$D^{\text{U}}$ [MW]},
        legend style={
        color=white,
        draw=white,
        fill=none,
        only marks,
        legend cell align={left},
        font=\footnotesize,
        legend pos=north west},
        width=8cm,
        height=4.5cm]
\addplot [color_ConvD, mark=diamond, line width=0.3mm, mark size=4] table [x index = 0, y index = 1] {plot_data/req_probabilistic_mark.dat};
\addlegendentry{Probabilistic}
\addplot [color_StochD, mark=square, line width=0.3mm, mark size=3.5] table [x index = 0, y index = 1] {plot_data/req_stochastic_mark.dat};
\addlegendentry{Stochastic}
\addplot [color_bilevel, mark=triangle, line width=0.3mm, mark size=4] table [x index = 0, y index = 1] {plot_data/req_optimal_mark.dat};
\addlegendentry{Enhanced}
\addplot3[surf,shader=interp] table [row sep=newline] {plot_data/surface_reserve.dat};
\addplot [white, line width=0.1mm, dashed, smooth] table [x index = 0, y index = 1] {plot_data/req_stochastic_xline.dat};
\addplot [white, line width=0.1mm, dashed, smooth] table [x index = 0, y index = 1] {plot_data/req_stochastic_yline.dat};
\addplot [white, line width=0.1mm, dashed, smooth] table [x index = 0, y index = 1] {plot_data/req_optimal_xline.dat};
\addplot [white, line width=0.1mm, dashed, smooth] table [x index = 0, y index = 1] {plot_data/req_optimal_yline.dat};
\addplot [white, line width=0.1mm, dashed, smooth] table [x index = 0, y index = 1] {plot_data/req_probabilistic_xline.dat};
\addplot [white, line width=0.1mm, dashed, smooth] table [x index = 0, y index = 1] {plot_data/req_probabilistic_yline.dat};
\end{axis}
\end{tikzpicture}
}
        \label{fig:ex1}
    }
    \subfigure[Day-ahead energy cost]
    {
\resizebox{8cm}{!}{%
\begin{tikzpicture}[thick,scale=1]
\begin{axis}[view={0}{90},
colormap/viridis,
        xmin=50,
        xmax=350,
        ymax=160,
        ymin=30,
        colorbar,
        colorbar style=     {
        yticklabel style=   {
        try min ticks=6,
        /pgf/number format/.cd,
            fixed,
            fixed zerofill,
            precision=1,
        /tikz/.cd
                            },
        ylabel={Cost [\$1000]},
        y label style=  {
        rotate=0,
        yshift = -80,
                        },
        ymin=21.9548400,
        ymax=23.0116600,
                            },
        ylabel={$D^{\text{D}}$ [MW]},
        xlabel={$D^{\text{U}}$ [MW]},
        legend style={
        draw=gray,
        fill=none,
        only marks,
        legend cell align={left},
        font=\footnotesize,
        legend pos=north west},
        width=8cm,
        height=4.5cm]
\addplot [color_bilevel, mark=triangle, line width=0.3mm, mark size=4] table [x index = 0, y index = 1] {plot_data/req_optimal_mark.dat};
\addplot [color_ConvD, mark=diamond, line width=0.3mm, mark size=4] table [x index = 0, y index = 1] {plot_data/req_probabilistic_mark.dat};
\addplot [color_StochD!70, mark=square, line width=0.3mm, mark size=3.5] table [x index = 0, y index = 1] {plot_data/req_stochastic_mark.dat};
\addplot3[surf,shader=interp] table [row sep=newline] {plot_data/surface_dayahead.dat};
\addplot [white, line width=0.1mm, dashed, smooth] table [x index = 0, y index = 1] {plot_data/req_stochastic_xline.dat};
\addplot [white, line width=0.1mm, dashed, smooth] table [x index = 0, y index = 1] {plot_data/req_stochastic_yline.dat};
\addplot [white, line width=0.1mm, dashed, smooth] table [x index = 0, y index = 1] {plot_data/req_optimal_xline.dat};
\addplot [white, line width=0.1mm, dashed, smooth] table [x index = 0, y index = 1] {plot_data/req_optimal_yline.dat};
\addplot [white, line width=0.1mm, dashed, smooth] table [x index = 0, y index = 1] {plot_data/req_probabilistic_xline.dat};
\addplot [white, line width=0.1mm, dashed, smooth] table [x index = 0, y index = 1] {plot_data/req_probabilistic_yline.dat};
\end{axis}
\end{tikzpicture}
}
        \label{fig:ex2}
    }
    \\
    \subfigure[Expected re-dispatch cost]
    {
\resizebox{8cm}{!}{%
\begin{tikzpicture}[thick,scale=1]
\begin{axis}[view={0}{90},
colormap/viridis,
        xmin=50,
        xmax=350,
        ymax=160,
        ymin=30,
        colorbar,
        colorbar style=     {
        yticklabel style=   {
        try min ticks=6,
        /pgf/number format/.cd,
            fixed,
            fixed zerofill,
            precision=1,
        /tikz/.cd
                            },
        ylabel={Cost [\$1000]},
        y label style=  {
        rotate=0,
        yshift = -80,
                        },
        ymin=0.0392900,
        ymax=7.8859200,
                            },
        ylabel={$D^{\text{D}}$ [MW]},
        xlabel={$D^{\text{U}}$ [MW]},
        legend style={
        draw=gray,
        fill=none,
        only marks,
        legend cell align={left},
        font=\footnotesize,
        legend pos=north west},
        width=8cm,
        height=4.5cm]
\addplot [color_bilevel, mark=triangle, line width=0.3mm, mark size=4] table [x index = 0, y index = 1] {plot_data/req_optimal_mark.dat};
\addplot [color_ConvD, mark=diamond, line width=0.3mm, mark size=4] table [x index = 0, y index = 1] {plot_data/req_probabilistic_mark.dat};
\addplot [color_StochD, mark=square, line width=0.3mm, mark size=3.5] table [x index = 0, y index = 1] {plot_data/req_stochastic_mark.dat};
\addplot3[surf,shader=interp] table [row sep=newline] {plot_data/surface_realtime.dat};
\addplot [white, line width=0.1mm, dashed, smooth] table [x index = 0, y index = 1] {plot_data/req_stochastic_xline.dat};
\addplot [white, line width=0.1mm, dashed, smooth] table [x index = 0, y index = 1] {plot_data/req_stochastic_yline.dat};
\addplot [white, line width=0.1mm, dashed, smooth] table [x index = 0, y index = 1] {plot_data/req_optimal_xline.dat};
\addplot [white, line width=0.1mm, dashed, smooth] table [x index = 0, y index = 1] {plot_data/req_optimal_yline.dat};
\addplot [white, line width=0.1mm, dashed, smooth] table [x index = 0, y index = 1] {plot_data/req_probabilistic_xline.dat};
\addplot [white, line width=0.1mm, dashed, smooth] table [x index = 0, y index = 1] {plot_data/req_probabilistic_yline.dat};
\end{axis}
\end{tikzpicture}
}
        \label{fig:ex3}
    }
    \subfigure[Expected total system cost]
    {
\resizebox{8cm}{!}{%
\begin{tikzpicture}[thick,scale=1]
\begin{axis}[view={0}{90},
colormap/viridis,
        xmin=50,
        xmax=350,
        ymax=160,
        ymin=30,
        colorbar,
        colorbar style=     {
        yticklabel style=   {
        try min ticks=6,
        /pgf/number format/.cd,
            fixed,
            fixed zerofill,
            precision=1,
        /tikz/.cd
                            },
        ylabel={Cost [\$1000]},
        y label style=  {
        rotate=0,
        yshift = -80,
                        },
        ymin=24.4,
        ymax=30,
                            },
        ylabel={$D^{\text{D}}$ [MW]},
        xlabel={$D^{\text{U}}$ [MW]},
        legend style={
        draw=gray,
        fill=none,
        only marks,
        legend cell align={left},
        font=\footnotesize,
        legend pos=north west},
        width=8cm,
        height=4.5cm]
\addplot [color_bilevel, mark=triangle, line width=0.3mm, mark size=4] table [x index = 0, y index = 1] {plot_data/req_optimal_mark.dat};
\addplot [color_ConvD, mark=diamond, line width=0.3mm, mark size=4] table [x index = 0, y index = 1] {plot_data/req_probabilistic_mark.dat};
\addplot [color_StochD, mark=square, line width=0.3mm, mark size=3.5] table [x index = 0, y index = 1] {plot_data/req_stochastic_mark.dat};
\addplot3[surf,shader=interp] table [row sep=newline] {plot_data/surface_total.dat};
\addplot [white, line width=0.1mm, dashed, smooth] table [x index = 0, y index = 1] {plot_data/req_stochastic_xline.dat};
\addplot [white, line width=0.1mm, dashed, smooth] table [x index = 0, y index = 1] {plot_data/req_stochastic_yline.dat};
\addplot [white, line width=0.1mm, dashed, smooth] table [x index = 0, y index = 1] {plot_data/req_optimal_xline.dat};
\addplot [white, line width=0.1mm, dashed, smooth] table [x index = 0, y index = 1] {plot_data/req_optimal_yline.dat};
\addplot [white, line width=0.1mm, dashed, smooth] table [x index = 0, y index = 1] {plot_data/req_probabilistic_xline.dat};
\addplot [white, line width=0.1mm, dashed, smooth] table [x index = 0, y index = 1] {plot_data/req_probabilistic_yline.dat};
\end{axis}
\end{tikzpicture}
}
        \label{fig:ex4}
    }
    \caption{Impact of downward $D^{\text{D}}$ and upward $D^{\text{U}}$ operating reserve requirements on the reserve (a), day-ahead (b), expected re-dispatch (c) and expected total (d) costs in the conventional procedure (\ref{prob:reserve_clearing})-(\ref{prob:real_time_clearing}). The color density indicates the cost at the considered trading floor.}
    \label{fig:cost_comparison}
\squeezeup
\end{figure*}
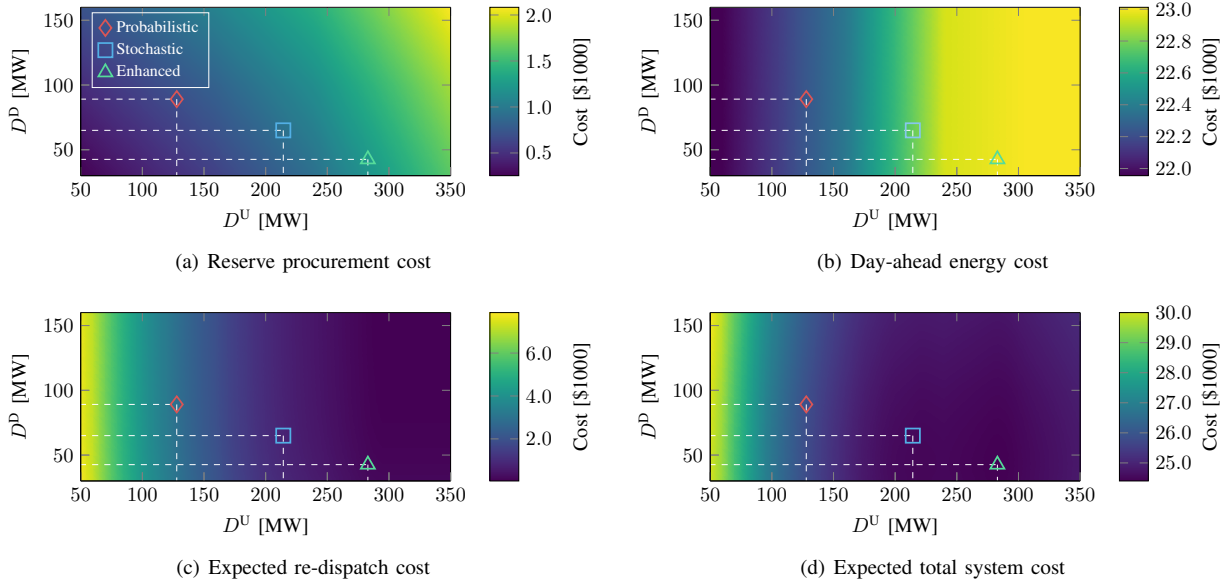

\subsection{Description of the test system} \label{Desription_of_test_system}
To assess the performance of the different reserve determination models, a modified version of the IEEE 24-Bus RTS \cite{Ordoudis_2016} is employed. The system consists of 34 transmission lines, 17 loads and 12 conventional generation units. The total generation capacity amounts to 3,375 MW, from which 1,100 MW is flexible generation that can provide upward and downward reserves.  We set upward reserve capacity price offers to be 30\% of the marginal costs. Price offers for downward reserve capacity price offers are selected such that they compensate for the potential financial deficit induced by a loss-making position in the day-ahead market. \textcolor{blue}{We should note that this is only a heuristic approach to address the possibility that some flexible producers incur financial losses due to their combined positions in the reserve capacity and day-ahead energy markets. This situation may emerge if the downward reserve capacity $R_i^{D*}$ awarded to a generator, and in turn imposed as a lower bound in the day-ahead market constraint (2c), forces this unit to produce even if the day-ahead energy price is lower than its marginal production cost. This pitfall results from the separation of reserve capacity and energy markets in the European framework. In turn, the physical coupling of these two products is accounted for internally in the trading strategies of the market participants when they submit their price-quantity offers in the corresponding markets according to their risk appetite. A detailed study of this issue constitutes a separate research topic and lies out of the scope of this work, but the interested reader is referred to \cite{SWIDER20071297} and \cite{1525135} for further information on this topic.} Apart from conventional generators, there are six wind farms bidding at zero marginal cost and sited as explained in \cite{Ordoudis_2016}. We consider a 24-hour load profile with a peak value of 2,650 MW obtained from \cite{Ordoudis_2016}. The loads are assumed to be inelastic with the value of lost load equal to \$500/MW for all operating hours. The relevant GAMS codes and simulation data are provided in the electronic companion of the paper \cite{companion}.

\textcolor{blue}{All simulations are carried out using a standard PC with Intel Core i5 CPU with a clock rate of 2.7 GHz requiring no more than 8GB of RAM. The CPU time required to solve the conventional model \eqref{prob:reserve_clearing}-\eqref{prob:real_time_clearing}, stochastic model \eqref{stochastic_dis} and bilevel model \eqref{prob:bilevel_clearing} in Sections \ref{Impact_of_reserve_requirements}--\ref{Optimal_zonal_reserve} is kept below 30s when solving per operating hour. The sequential market with unit commitment and inter-temporal constraints is solved in less than a minute in Section \ref{24RTS_with_UC_conctrains}. The CPU time corresponding to the last case study is reported separately in Section \ref{IEEE96study}.
}

\subsection{Impact of reserve requirements on expected system cost} \label{Impact_of_reserve_requirements}
\textcolor{blue}{In this section we assess the expected cost of operating the power system under the conventional market setup \eqref{prob:reserve_clearing}-\eqref{prob:real_time_clearing}, when this is fed with the reserve requirements determined by different approaches for reserve dimensioning, including our proposal. To this end, we consider the time period corresponding to the peak-load hour. Besides, the capacity of each wind power farm is set to 100 MW. Next we discuss the results linked to each reserve dimensioning approach:}
\begin{enumerate}
    \item The \textit{probabilistic approach} \textcolor{ForestGreen}{defines the reserve} requirements from the \textcolor{WildStrawberry}{predictive} cumulative distribution function (CDF) $F$ of the total wind power portfolio, as the distance between the expected wind power production $\widehat{W}$ and a specified quantile $q^{(\alpha)} = F^{-1}(\alpha)$ \textcolor{WildStrawberry}{with nominal proportion} $\alpha \in [0,1]$. \textcolor{blue}{This approach resembles the state-of-the-art reserve-dimensioning processes employed by European system operators using probabilistic forecast information \cite{6299425}}. \textcolor{WildStrawberry}{For a reliability level $\xi = \overline{\alpha} - \underline{\alpha}$}, the upward and downward reserve needs are dimensioned as follows:
    \begin{subequations}
    \begin{align}
        & D^{\text{U}} = \widehat{W} - F^{-1}(\overline{\alpha}), \\
        & D^{\text{D}} = F^{-1}(\underline{\alpha}) - \widehat{W}.
    \end{align}
    \end{subequations}
    We initially consider \textcolor{WildStrawberry}{$\underline{\alpha} = 5\%$ and $\overline{\alpha} = 1 - \underline{\alpha} = 95\%$ corresponding to a reliability level $\xi=$ 90\%.} The resulting  requirements amount to 127.9 MW and 89.1 MW for upward and downward reserves, respectively.
    \item \textit{The stochastic approach} \textcolor{ForestGreen}{derives the reserve} requirements from the stochastic dispatch model (\ref{stochastic_dis}). These requirements are equal to 214.3 MW for upward and 65.0 MW for downward reserves, respectively.
    \item \textit{The enhanced approach} computes the \textcolor{ForestGreen}{reserve} requirements using the proposed reserve determination model (\ref{prob:bilevel_clearing}). Resulting reserve needs amount to 282.9 MW and 42.6 MW for upward and downward reserves, respectively.
\end{enumerate}

The expected total system costs resulting from the implementation of the \textcolor{ForestGreen}{probabilistic, stochastic and enhanced operating reserve} approaches are \$25,890, \$24,531 and \$24,408, respectively. The total cost break-down is shown in Fig. \ref{fig:cost_comparison}, which demonstrates the impact of the reserve requirements on the cost of the different trading floors in the conventional dispatch procedure. Figure \ref{fig:ex1} shows that the reserve needs computed using the proposed model result in the highest reserve procurement cost among the different approaches, mainly due to a larger volume of upward reserve provision. In turn, efficient flexible generation that could be scheduled in the day-ahead market is now set aside to provide upward reserves. Considering that the price offers for upward reserve are proportional to the day-ahead price offers, the withdrawal of these resources increases the day-ahead energy cost, as shown in Fig. \ref{fig:ex2}. Nonetheless, the benefits of the enhanced approach realize in real-time operation as the re-dispatch cost is lower compared to that yielded by the probabilistic and stochastic approaches as illustrated in Fig. \ref{fig:ex3}. As a result, the minimum of the expected total costs is achieved with the enhanced approach as demonstrated by Fig. \ref{fig:ex4}.

\textcolor{WildStrawberry}{Increasing the reliability level $\xi$ in the the probabilistic approach may have a positive impact on the performance of the conventional model.} However, Table \ref{cost_break_down} shows that this approach never yields the expected cost provided by the proposed model, since the probabilistic approach sets the requirements disregarding their impact on the subsequent operations, including potential wind spillage and load shedding.
\textcolor{WildStrawberry}{On the contrary}, the proposed model finds \textcolor{WildStrawberry}{the optimal} trade-off between reserve procurement and real-time re-dispatch decisions that minimizes the total expected  system cost. In this particular case, our model allows more wind curtailment to reduce downward reserve procurement cost.

\textcolor{WildStrawberry}{Regarding the stochastic model, it should be noted that even though} reserve requirements are set anticipating the real-time cost, reserve procurement and day-ahead energy schedules are obtained by a co-optimization of these products \textcolor{WildStrawberry}{that is incompatible with the European} market structure. As a result, the requirements provided by the stochastic approach lead to larger amounts of load shedding, highlighting that they are \textcolor{WildStrawberry}{practically} sub-optimal in a sequential dispatch procedure.

\begin{table}[]
\centering
\caption{Cost break-down resulting from the implementation of a range of probabilistic requirements and enhanced requirements.}
\label{cost_break_down}
\tabcolsep=0.06cm
\resizebox{0.48\textwidth}{!}{\begin{tabular}{l|ccccc|c}
\specialrule{1pt}{1pt}{1pt}
\multicolumn{1}{c|}{\multirow{3}{*}{Approach}} & \multicolumn{5}{c|}{Probabilistic approach} & \multirow{3}{*}{\begin{tabular}[c]{@{}c@{}}Enhanced\\ approach\end{tabular}} \\
\multicolumn{1}{c|}{} & \multicolumn{5}{c|}{ \textcolor{WildStrawberry}{Quantiles $q^{(\underline{\alpha},\overline{\alpha})}$} of wind CDF} &  \\
\cline{2-6}
\multicolumn{1}{c|}{} & $q^{(05/95)}$ & $q^{(04/96)}$ & $q^{(03/97)}$ & $q^{(02/98)}$ & $q^{(01/99)}$ &  \\
\specialrule{1pt}{1pt}{1pt}
Requirements $D^{\text{U/D}}$ {[}MW{]} & 128/89 & 168/91 & 205/93 & 210/94 & 283/169 & 283/43 \\
Exp. total cost {[}\$1000{]} & 25.89 & 24.99 & 24.62 & 24.61 & 24.78 & 24.40 \\
-- \textit{Reserve} & 0.69 & 0.84 & 0.99 & 1.01 & 1.70 & 1.24 \\
-- \textit{Day-ahead} & 22.24 & 22.43 & 22.70 & 22.74 & 22.99 & 22.99 \\
-- \textit{Real-time} & 2.96 & 1.72 & 0.93 & 0.86 & 0.88 & 0.18 \\
\specialrule{1pt}{1pt}{1pt}
\end{tabular}}

\end{table}

\subsection{Approximating the stochastic dispatch solution} \label{Approximating_the_stoch}

\begin{figure}[]
\center
\begin{tikzpicture}[thick,scale=1]
\pgfplotsset{ymin=285, ymax=500, xmin=0, xmax=50, try min ticks=5}
\begin{axis}[xlabel near ticks, ylabel near ticks, ylabel={Expected cost [\$1000]}, xlabel={Wind penetration [\% of peak load]}, label style={font=\footnotesize},
 tick label style={font=\scriptsize}, legend pos = south west, legend style={draw=none, name = ConvD_leg_pos, font=\scriptsize, legend cell align={left}}, width=8.6cm, height=4.5cm]
    \addplot [color_ConvD!120, line width=0.15mm, mark=*, mark size=1.1, smooth, name path=A] table [x index = 0, y index = 1] {plot_data/Cost_data_wind_share.dat};
    \addlegendentry{Probabilistic solution}
    \addplot [color_bilevel!120, line width=0.15mm, mark=diamond*, mark size=1.25, smooth] table [x index = 0, y index = 3] {plot_data/Cost_data_wind_share.dat};
    \addlegendentry{Enhanced solution}
    \addplot [color_StochD!120, line width=0.15mm, mark=square*, mark size=1.1, smooth] table [x index = 0, y index = 4] {plot_data/Cost_data_wind_share.dat};
    \addlegendentry{Ideal solution}
    \addplot [color_ConvD!120, line width=0.15mm, mark=*, mark size=1.1, smooth, name path=B] table [x index = 0, y index = 2] {plot_data/Cost_data_wind_share.dat};
    \addplot[pattern=crosshatch dots, pattern color=color_ConvD!50]fill between[of=A and B];
\end{axis}
\end{tikzpicture}
\caption{Expected daily operating cost as a function of wind penetration.}
\label{costs_diff_wind}
\end{figure}
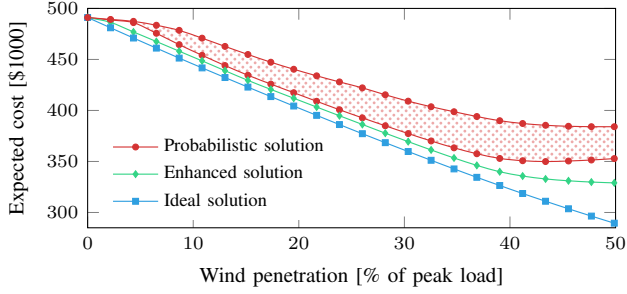

\begin{figure}
    \centering
    \subfigure[Upward reserve procurement]
    {
\begin{tikzpicture}[thick,scale=1]
  \pgfplotsset{ybar stacked, ymin=0, ymax=600, xmin=0.5, xmax=5.5, xtick=data}
  \begin{axis}[xlabel near ticks,ylabel near ticks,bar shift=-10pt, xlabel={Wind penetration [\% of peak demand]}, ylabel={$R_{i}^{\text{U}}$ [MW]}, label style={font=\footnotesize}, xticklabels={10.84,21.69,30.36,41.2,49.88}, tick label style={font=\footnotesize}, legend pos = north west, legend style={draw=none,  font=\scriptsize,name = ConvD_leg_pos}, width=8cm, height=4cm]
    \addplot [white!120, line width=0.05mm, fill=color_ConvD!30] table [x index = 0, y index = 1] {plot_data/Prob_UR.dat};
    \addplot [white!120, line width=0.05mm, fill=color_ConvD!40] table [x index = 0, y index = 2] {plot_data/Prob_UR.dat};
    \addplot [white!120, line width=0.05mm, fill=color_ConvD!50] table [x index = 0, y index = 3] {plot_data/Prob_UR.dat};
    \addplot [white!120, line width=0.05mm, fill=color_ConvD!60] table [x index = 0, y index = 4] {plot_data/Prob_UR.dat};
    \addplot [white!120, line width=0.05mm, fill=color_ConvD!70] table [x index = 0, y index = 5] {plot_data/Prob_UR.dat};
    \addplot [white!120, line width=0.05mm, fill=color_ConvD!80] table [x index = 0, y index = 6] {plot_data/Prob_UR.dat};
    \addplot [white!120, line width=0.05mm, fill=color_ConvD!110] table [x index = 0, y index = 7] {plot_data/Prob_UR.dat};
    \addplot [white!120, line width=0.05mm, fill=color_ConvD!130] table [x index = 0, y index = 8] {plot_data/Prob_UR.dat};
    \addplot [white!120, line width=0.05mm, fill=color_ConvD!150] table [x index = 0, y index = 9] {plot_data/Prob_UR.dat};
    \draw [line width=0.05mm, draw=color_ConvD!150]
        (axis cs: 0.725,20) -- (axis cs: 0.725,525);
    \draw [line width=0.05mm, draw=color_ConvD!150]
        (axis cs: 0.725,525) -- (axis cs: 1.5,525)
        node[right] {\tiny Probabilistic solution, $\alpha = 0.05$ };
    \draw [line width=0.05mm, draw=color_bilevel!150]
        (axis cs: 1,20) -- (axis cs: 1,470);
    \draw [line width=0.05mm, draw=color_bilevel!150]
        (axis cs: 1,470) -- (axis cs: 1.5,470)
        node[right] {\tiny Enhanced solution };
    \draw [line width=0.05mm, draw=color_StochD!150]
        (axis cs: 1.275,20) -- (axis cs: 1.275,415);
    \draw [line width=0.05mm, draw=color_StochD!150]
        (axis cs: 1.275,415) -- (axis cs: 1.5,415)
        node[right] {\tiny Ideal solution };
  \end{axis}
  \begin{axis}[hide axis, xticklabels={}, yticklabels={},legend pos = north west , legend style={draw=none,  font=\scriptsize, name = Bilevel_leg_pos, at = {([yshift = 2mm]ConvD_leg_pos.south west)},
      anchor = north west}, width=8cm, height=4cm]
    \addplot [white!120, line width=0.05mm, fill=color_bilevel!30] table [x index = 0, y index = 1] {plot_data/Enhanced_UR.dat};
    \addplot [white!120, line width=0.05mm, fill=color_bilevel!40] table [x index = 0, y index = 2] {plot_data/Enhanced_UR.dat};
    \addplot [white!120, line width=0.05mm, fill=color_bilevel!50] table [x index = 0, y index = 3] {plot_data/Enhanced_UR.dat};
    \addplot [white!120, line width=0.05mm, fill=color_bilevel!60] table [x index = 0, y index = 4] {plot_data/Enhanced_UR.dat};
    \addplot [white!120, line width=0.05mm, fill=color_bilevel!70] table [x index = 0, y index = 5] {plot_data/Enhanced_UR.dat};
    \addplot [white!120, line width=0.05mm, fill=color_bilevel!80] table [x index = 0, y index = 6] {plot_data/Enhanced_UR.dat};
    \addplot [white!120, line width=0.05mm, fill=color_bilevel!110] table [x index = 0, y index = 7] {plot_data/Enhanced_UR.dat};
    \addplot [white!120, line width=0.05mm, fill=color_bilevel!130] table [x index = 0, y index = 8] {plot_data/Enhanced_UR.dat};
    \addplot [white!120, line width=0.05mm, fill=color_bilevel!150] table [x index = 0, y index = 9] {plot_data/Enhanced_UR.dat};
  \end{axis}
  \begin{axis}[hide axis, bar shift = 10pt, xticklabels={}, yticklabels={},legend pos = north west , legend style={draw=none,  font=\scriptsize, name = StochD_leg_pos, at = {([yshift = 2mm]Bilevel_leg_pos.south west)},
      anchor = north west}, width=8cm, height=4cm]
    \addplot [white!120, line width=0.05mm, fill=color_StochD!30] table [x index = 0, y index = 1] {plot_data/Indeal_UR.dat};
    \addplot [white!120, line width=0.05mm, fill=color_StochD!40] table [x index = 0, y index = 2] {plot_data/Indeal_UR.dat};
    \addplot [white!120, line width=0.05mm, fill=color_StochD!50] table [x index = 0, y index = 3] {plot_data/Indeal_UR.dat};
    \addplot [white!120, line width=0.05mm, fill=color_StochD!60] table [x index = 0, y index = 4] {plot_data/Indeal_UR.dat};
    \addplot [white!120, line width=0.05mm, fill=color_StochD!70] table [x index = 0, y index = 5] {plot_data/Indeal_UR.dat};
    \addplot [white!120, line width=0.05mm, fill=color_StochD!80] table [x index = 0, y index = 6] {plot_data/Indeal_UR.dat};
    \addplot [white!120, line width=0.05mm, fill=color_StochD!110] table [x index = 0, y index = 7] {plot_data/Indeal_UR.dat};
    \addplot [white!120, line width=0.05mm, fill=color_StochD!130] table [x index = 0, y index = 8] {plot_data/Indeal_UR.dat};
    \addplot [white!120, line width=0.05mm, fill=color_StochD!150] table [x index = 0, y index = 9] {plot_data/Indeal_UR.dat};
  \end{axis}
\end{tikzpicture}
    }
    \\
    \subfigure[Downward reserve procurement]
    {
\begin{tikzpicture}[thick,scale=1]
  \pgfplotsset{ybar stacked, ymin=0, ymax=220, xmin=0.5, xmax=5.5, xtick=data}
  \begin{axis}[xlabel near ticks,ylabel near ticks,bar shift=-10pt, xlabel={Wind penetration [\% of peak demand]}, ylabel={$R_{i}^{\text{D}}$ [MW]}, label style={font=\footnotesize}, xticklabels={10.84,21.69,30.36,41.2,49.88}, tick label style={font=\footnotesize}, legend pos = north west, legend style={draw=none,  font=\scriptsize,name = ConvD_leg_pos}, width=8cm, height=4cm]
    \addplot [white!120, line width=0.05mm, fill=color_ConvD!50] table [x index = 0, y index = 1] {plot_data/Prob_DR.dat};
    \addplot [white!120, line width=0.05mm, fill=color_ConvD!60] table [x index = 0, y index = 2] {plot_data/Prob_DR.dat};
    \addplot [white!120, line width=0.05mm, fill=color_ConvD!70] table [x index = 0, y index = 3] {plot_data/Prob_DR.dat};
    \addplot [white!120, line width=0.05mm, fill=color_ConvD!80] table [x index = 0, y index = 4] {plot_data/Prob_DR.dat};
    \addplot [white!120, line width=0.05mm, fill=color_ConvD!90] table [x index = 0, y index = 5] {plot_data/Prob_DR.dat};
    \addplot [white!120, line width=0.05mm, fill=color_ConvD!110] table [x index = 0, y index = 6] {plot_data/Prob_DR.dat};
    \addplot [white!120, line width=0.05mm, fill=color_ConvD!130] table [x index = 0, y index = 7] {plot_data/Prob_DR.dat};
    \addplot [white!120, line width=0.05mm, fill=color_ConvD!140] table [x index = 0, y index = 8] {plot_data/Prob_DR.dat};
    \addplot [white!120, line width=0.05mm, fill=color_ConvD!150] table [x index = 0, y index = 9] {plot_data/Prob_DR.dat};
  \end{axis}
  \begin{axis}[hide axis, xticklabels={}, yticklabels={},legend pos = north west , legend style={draw=none,  font=\scriptsize, name = Bilevel_leg_pos, at = {([yshift = 2mm]ConvD_leg_pos.south west)},
      anchor = north west}, width=8cm, height=4cm]
    \addplot [white!120, line width=0.05mm, fill=color_bilevel!50] table [x index = 0, y index = 1] {plot_data/Enhanced_DR.dat};
    \addplot [white!120, line width=0.05mm, fill=color_bilevel!60] table [x index = 0, y index = 2] {plot_data/Enhanced_DR.dat};
    \addplot [white!120, line width=0.05mm, fill=color_bilevel!70] table [x index = 0, y index = 3] {plot_data/Enhanced_DR.dat};
    \addplot [white!120, line width=0.05mm, fill=color_bilevel!80] table [x index = 0, y index = 4] {plot_data/Enhanced_DR.dat};
    \addplot [white!120, line width=0.05mm, fill=color_bilevel!90] table [x index = 0, y index = 5] {plot_data/Enhanced_DR.dat};
    \addplot [white!120, line width=0.05mm, fill=color_bilevel!110] table [x index = 0, y index = 6] {plot_data/Enhanced_DR.dat};
    \addplot [white!120, line width=0.05mm, fill=color_bilevel!130] table [x index = 0, y index = 7] {plot_data/Enhanced_DR.dat};
    \addplot [white!120, line width=0.05mm, fill=color_bilevel!140] table [x index = 0, y index = 8] {plot_data/Enhanced_DR.dat};
    \addplot [white!120, line width=0.05mm, fill=color_bilevel!150] table [x index = 0, y index = 9] {plot_data/Enhanced_DR.dat};
  \end{axis}
  \begin{axis}[hide axis, bar shift = 10pt, xticklabels={}, yticklabels={},legend pos = north west , legend style={draw=none,  font=\scriptsize, name = StochD_leg_pos, at = {([yshift = 2mm]Bilevel_leg_pos.south west)},
      anchor = north west}, width=8cm, height=4cm]
    \addplot [white!120, line width=0.05mm, fill=color_StochD!50] table [x index = 0, y index = 1] {plot_data/Indeal_DR.dat};
    \addplot [white!120, line width=0.05mm, fill=color_StochD!60] table [x index = 0, y index = 2] {plot_data/Indeal_DR.dat};
    \addplot [white!120, line width=0.05mm, fill=color_StochD!70] table [x index = 0, y index = 3] {plot_data/Indeal_DR.dat};
    \addplot [white!120, line width=0.05mm, fill=color_StochD!80] table [x index = 0, y index = 4] {plot_data/Indeal_DR.dat};
    \addplot [white!120, line width=0.05mm, fill=color_StochD!90] table [x index = 0, y index = 5] {plot_data/Indeal_DR.dat};
    \addplot [white!120, line width=0.05mm, fill=color_StochD!120] table [x index = 0, y index = 6] {plot_data/Indeal_DR.dat};
    \addplot [white!120, line width=0.05mm, fill=color_StochD!130] table [x index = 0, y index = 7] {plot_data/Indeal_DR.dat};
    \addplot [white!120, line width=0.05mm, fill=color_StochD!140] table [x index = 0, y index = 8] {plot_data/Indeal_DR.dat};
    \addplot [white!120, line width=0.05mm, fill=color_StochD!150] table [x index = 0, y index = 9] {plot_data/Indeal_DR.dat};
  \end{axis}
\end{tikzpicture}
    }
    \caption{Reserve procurement from nine flexible generating units for the peak-load hour and different wind penetration levels. Color density ranks generation units according to the reserve capacity price offers.}
    \label{Reserve_procur}
\squeezeup
\end{figure}
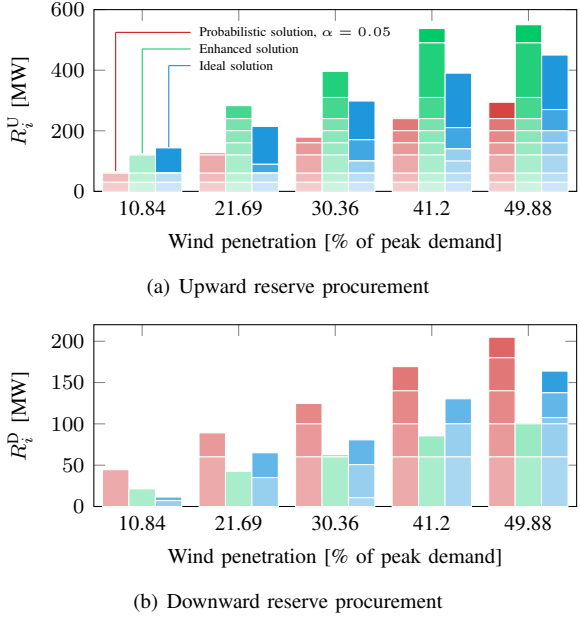

We now investigate to what extent the reserve requirements computed with the proposed model are capable of approximating the ideal stochastic solution within the sequential dispatch procedure. To this end, we compare expected daily system cost of three optimization models for different wind power penetration levels, defined as the ratio between the installed capacity of the entire wind power portfolio and the peak load. The first model represents the sequential market clearing \eqref{prob:reserve_clearing}-\eqref{prob:real_time_clearing} with reserve requirements computed with the probabilistic approach for a range of \textcolor{WildStrawberry}{ reliability levels $\xi \in [0.9,1]$}. The second model \textcolor{WildStrawberry}{also follows the} sequential market procedure with reserve requirements computed with the proposed model \eqref{prob:bilevel_clearing}. The third one is the stochastic ideal dispatch model \eqref{stochastic_dis} that theoretically attains maximum cost-efficiency, and therefore it is used as a lower bound of the expected system cost. \textcolor{blue}{It is worth noting the different role that the stochastic dispatch model plays in this part of the case study, compared to the previous Section \ref{Impact_of_reserve_requirements}. Here, we assume that the solution of the stochastic dispatch model will be implemented as the actual system schedule, presuming that the conventional market setup is replaced with its ideal stochastic counterpart. This is different from the application of the stochastic dispatch model \eqref{stochastic_dis} as a reserve-dimensioning approach in Section \ref{Impact_of_reserve_requirements}, where we considered that all trading floors are settled according to the prevailing European market model.}

Figure \ref{costs_diff_wind} depicts the daily operating cost as a function of the wind power penetration level for the three models. The setting of the reserve requirements provided by the proposed model always results in a lower expected cost than the implementation of the requirements under the probabilistic approach. This figure further indicates that these reserve requirements \textcolor{WildStrawberry}{efficiently approximate the stochastic ideal solution even for a high penetration of wind power.}

Figure \ref{Reserve_procur} provides further insights on the difference between the solutions of the three models. Particularly, it shows the procurement of upward and downward reserves from specific flexible units ranked according to their reserve capacity price offers, i.e., from cheap to more expensive units distinguished by increasing color densities. The proposed model controls the trade-off between reserve and real-time costs, ensuring  adequate upward reserves to minimize the amount of load shedding and enough downward reserves to prevent wind spillage. In contrast, the probabilistic approach underestimates upward reserve needs, while it overestimates downward reserve requirements.

The enhanced solution for the reserve requirements deviates significantly from the ideal solution given that the stochastic model has more degrees of freedom, i.e., it controls not only the sufficiency of the reserve requirements but also their allocation among the flexible generators. This results in reserve procurement  being `generator-specific' which prevents network congestion within the reserve control area. In attempt to minimize expensive balancing actions, the stochastic model \textcolor{WildStrawberry}{may allocate reserves to more expensive units over cheaper providers}, violating the least-cost merit-order principle that is inherent in the conventional market design. \textcolor{ForestGreen}{As a consequence, the requirement imposed in our enhanced approach to respect the merit-order principle in the reserve capacity and day-ahead markets restricts the degree of approximation of the stochastic solution.}

\subsection{Optimal zonal reserve requirements allocation} \label{Optimal_zonal_reserve}
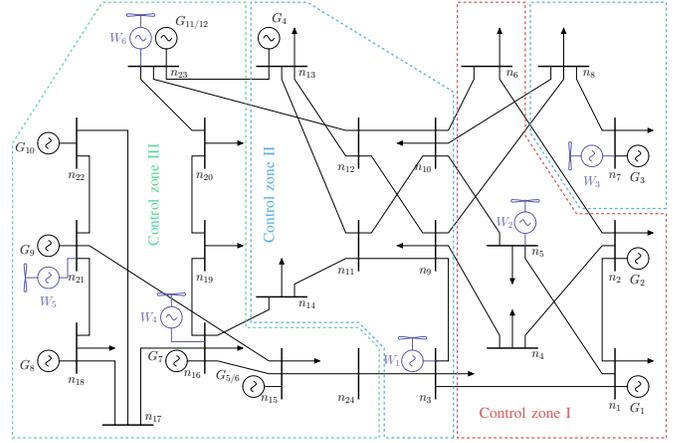
\begin{figure}
\centering
\resizebox{8.75cm}{!}{%
\begin{turn}{90}
\begin{circuitikz}[>=triangle 45]
\draw [ultra thick] (2,6) node[anchor=north, rotate = 270]{\Large$n_{1}$} -- (4,6);
\draw [ultra thick] (7,6) node[anchor=north, rotate = 270]{\Large$n_{2}$} -- (9,6);
\draw [ultra thick] (2,13) node[anchor=east, rotate = 270]{\Large$n_{3}$} -- (4,13);
\draw [ultra thick] (4,9) node[anchor=north, rotate = 270]{\Large$n_{4}$} -- (4,11);
\draw [ultra thick] (8,9) node[anchor=north, rotate = 270]{\Large$n_{5}$} -- (8,11);
\draw [ultra thick] (15,10) node[anchor=north, rotate = 270]{\Large$n_{6}$} -- (15,12);
\draw [ultra thick] (11,6) node[anchor=north, rotate = 270]{\Large$n_{7}$} -- (13,6);
\draw [ultra thick] (15,7) node[anchor=north, rotate = 270]{\Large$n_{8}$} -- (15,9);
\draw [ultra thick] (7,13) node[anchor=east, rotate = 270]{\Large$n_{9}$} -- (9,13);
\draw [ultra thick] (11,13) node[anchor=east, rotate = 270]{\Large$n_{10}$} -- (13,13);
\draw [ultra thick] (7,16) node[anchor=east, rotate = 270]{\Large$n_{11}$} -- (9,16);
\draw [ultra thick] (11,16) node[anchor=east, rotate = 270]{\Large$n_{12}$} -- (13,16);
\draw [ultra thick] (15,18) node[anchor=north, rotate = 270]{\Large$n_{13}$} -- (15,20);
\draw [ultra thick] (6,18) node[anchor=north, rotate = 270]{\Large$n_{14}$} -- (6,20);
\draw [ultra thick] (2,19) node[anchor=east, rotate = 270]{\Large$n_{15}$} -- (4,19);
\draw [ultra thick] (3,22) node[anchor=east, rotate = 270]{\Large$n_{16}$} -- (5,22);
\draw [ultra thick] (1,24) node[anchor=south, rotate = 270]{\Large$n_{17}$} -- (1,26);
\draw [ultra thick] (3,27) node[anchor=north, rotate = 270]{\Large$n_{18}$} -- (5,27);
\draw [ultra thick] (7,22) node[anchor=north, rotate = 270]{\Large$n_{19}$} -- (9,22);
\draw [ultra thick] (11,22) node[anchor=north, rotate = 270]{\Large$n_{20}$} -- (13,22);
\draw [ultra thick] (7,27) node[anchor=north, rotate = 270]{\Large$n_{21}$} -- (9,27);
\draw [ultra thick] (11,27) node[anchor=north, rotate = 270]{\Large$n_{22}$} -- (13,27);
\draw [ultra thick] (15,23) node[anchor=north, rotate = 270]{\Large$n_{23}$} -- (15,25);
\draw [ultra thick] (2,16) node[anchor=east, rotate = 270]{\Large$n_{24}$} -- (4,16);
\draw[->] (3.5,6) -- +(0,-1.5); 
\draw[->] (8.5,6) -- +(0,-1.5); 
\draw[->] (3,13) -- +(0,-1.5); 
\draw[->] (4,10) -- +(1.5,0); 
\draw[->] (8,10) -- +(-1.5,0); 
\draw[->] (15,10.5) -- +(1.5,0); 
\draw[->] (12.5,6) -- +(0,-1.5); 
\draw[->] (15,8) -- +(1.5,0); 
\draw[->] (8,13) -- +(0,1.5); 
\draw[->] (12,13) -- +(0,1.5); 
\draw[->] (15,18.5) -- +(1.5,0); 
\draw[->] (6,19) -- +(1.5,0); 
\draw[->] (4,22) -- +(0,-1.5); 
\draw[->] (3.5,19) -- +(0,-1.5); 
\draw[->] (4,27) -- +(0,-1.5); 
\draw[->] (8,22) -- +(0,-1.5); 
\draw[->] (12,22) -- +(0,-1.5); 
\draw (12.5,6) -- (12.5,6.5) -- (14.5,7.5) -- (15,7.5);  
\draw (15,10.5) -- (14.5,10.5) -- (8.5,6.5) -- (8.5,6);  
\draw (15,11.5) -- (14.5,11.5) -- (12.5,12.5) -- (12.5,13);  
\draw (8,10.5) -- (8.5,10.5) -- (11.5,12.5) -- (11.5,13);  
\draw (15,8.5) -- (14.5,8.5) -- (12,12.5) -- (12,13);  
\draw (15,8) -- (14.5,8) -- (8.5,12.5) -- (8.5,13);  
\draw (4,10.5) -- (4.5,10.5) -- (8,12.5) -- (8,13);  
\draw (4,9.5) -- (4.5,9.5) -- (8,6.5) -- (8,6);  
\draw (3.5,6) -- (3.5,6.5) -- (7.5,6.5) -- (7.5,6);  
\draw (3,6) -- (3,6.5) -- (7.5,9.5) -- (8,9.5);  
\draw (2.5,13) -- (2.5,12.5) -- (2.5,6.5) -- (2.5,6);  
\draw (3.5,13) -- (3.5,12.5) -- (7.5,12.5) -- (7.5,13);  
\draw (8.5,13) -- (8.5,13.5) -- (11.5,15.5) -- (11.5,16);  
\draw (11.5,13) -- (11.5,13.5) -- (8.5,15.5) -- (8.5,16);  
\draw (12.5,13) -- (12.5,16);  
\draw (7.5,13) -- (7.5,16);  
\draw (3,13) -- (3,16);  
\draw (3,19) -- (3,16);  
\draw (11.5,16) -- (11.5,16.5) -- (14.5,18.5) -- (15,18.5);  
\draw (8.5,16) -- (8.5,16.5) -- (14.5,19) -- (15,19);  
\draw (15,23.5) -- (14.5,23.5) -- (14.5,19.5) -- (15,19.5);  
\draw (12.5,16) -- (12.5,16.5) -- (14.5,24) -- (15,24);  
\draw (12.5,22) -- (12.5,22.5) -- (14.5,24.5) -- (15,24.5);  
\draw (6,18.5) -- (6.5,18.5) -- (7.5,16.5) -- (7.5,16);  
\draw (4.5,22) -- (4.5,21.5) -- (5.5,19.5) -- (6,19.5);  
\draw (3.5,22) -- (3.5,21.5) -- (3,19.5) -- (3,19); 
\draw (8,27) -- (8,26.5) -- (3.5,19.5) -- (3.5,19); 
\draw (1,25) -- (12.5,25) -- (12.5,27);  
\draw (8.5,22) -- (8.5,22.5) -- (11.5,22.5) -- (11.5,22);  
\draw (4.5,22) -- (4.5,22.5) -- (7.5,22.5) -- (7.5,22);  
\draw (1,24.5) -- (4,24.5) -- (4,22);  
\draw (1,25.5) -- (3.5,25.5) -- (3.5,27);  
\draw (4.5,27) -- (4.5,26.5) -- (7.5,26.5) -- (7.5,27);  
\draw (8.5,27) -- (8.5,26.5) -- (11.5,26.5) -- (11.5,27); 
\draw (2.5,5.5)  to [sV] (2.5,4.7) node[anchor=north east, rotate=270, yshift = -15]{\Large$G_1$};
\draw (2.5,5.5)--(2.5,6);
\draw (7.5,5.5)  to [sV] (7.5,4.7) node[anchor=north east, rotate=270, yshift = -15]{\Large$G_2$};
\draw (7.5,5.5)--(7.5,6);
\draw (11.5,5.5)  to [sV] (11.5,4.7) node[anchor=north east, rotate=270, yshift = -15]{\Large$G_3$};
\draw (11.5,5.5)--(11.5,6);
\draw (16.5,19.5) node[anchor= west, rotate=270, yshift = 10]{\Large $G_4$}  to [sV] (15.7,19.5) ;
\draw (15,19.5)--(15.7,19.5);
\draw (2.5,20.5) node[anchor=south east, rotate=270]{\Large$G_{5/6}$} to [sV] (2.5,19.7) ;
\draw (2.5,19)--(2.5,19.7);
\draw (3.5,23.5) node[anchor=south east, rotate=270, yshift = -5]{\Large$G_7$} to [sV] (3.5,22.7) ;
\draw (3.5,22)--(3.5,22.7);
\draw (3.5,28.5) node[anchor=south east, rotate=270, yshift = -15]{\Large$G_8$} to [sV] (3.5,27.7) ;
\draw (3.5,27)--(3.5,27.7);
\draw (8,28.5) node[anchor=south east, rotate=270, yshift = -15]{\Large$G_9$} to [sV] (8,27.7);
\draw (8,27)--(8,27.7);
\draw (12,28.5) node[anchor=south east, rotate=270, yshift = -15]{\Large$G_{10}$} to [sV] (12,27.7);
\draw (12,27)--(12,27.7);
\draw (16.5,23.5) node[anchor= west]{}  to [sV] (15.7,23.5);
\node[rotate=270] at (16.7,22.5) {\Large$G_{11/12}$};
\draw (15,23.5)--(15.7,23.5);
\draw[color = color_dark_blue] (3.5,13.5)   to [sV] (3.5,14.3); 
\draw[color = color_dark_blue] (3.5,13)--(3.5,13.5);
\draw[color = color_dark_blue] plot [smooth cycle] coordinates {(4.45,14.4) (4.4,13.9) (4.35,13.4) (4.45,13.4) (4.4,13.9) (4.35,14.4)};
\draw[color = color_dark_blue] (4.4,13.9)--(3.9,13.9);
\node[rotate=270] at (3.5,14.7) {\Large$\textcolor{color_dark_blue}{W_1}$};
\draw[color = color_dark_blue] (8.5,9.5) to [sV] (9.3,9.5); 
\draw[color = color_dark_blue] (8.5,9.5)--(8.0,9.5);
\draw[color = color_dark_blue] (9.3,9.5)--(9.8,9.5);
\draw[color = color_dark_blue] plot [smooth cycle] coordinates {(9.8,10) (9.75,9.5) (9.7,9.0) (9.8,9.0) (9.75,9.5) (9.7,10)};
\node[rotate=270] at (8.9,10.3) {\Large$\textcolor{color_dark_blue}{W_2}$};
\draw[color = color_dark_blue] (11.5,6.5) to [sV] (11.5,7.3); 
\draw[color = color_dark_blue] (11.5,6.0)--(11.5,6.5);
\draw[color = color_dark_blue] plot [smooth cycle] coordinates {(11,7.8) (11.5,7.75) (12.0,7.7) (12.0,7.8) (11.5,7.75) (11.0,7.7)};
\draw[color = color_dark_blue] (11.5,7.3)--(11.5,7.75);
\node[rotate=270] at (10.5,6.9) {\Large$\textcolor{color_dark_blue}{W_3}$};
\draw[color = color_dark_blue] (4.8,23.3) to [sV] (5.6,23.3); 
\draw[color = color_dark_blue] (4.25,22) -- (4.25,23.3) -- (4.8,23.3);
\draw[color = color_dark_blue] (5.6,23.3) -- (6.1,23.3);
\draw[color = color_dark_blue] plot [smooth cycle] coordinates {(6.15, 23.8) (6.1,23.3) (6.05,22.8) (6.15,22.8) (6.1,23.3) (6.05,23.8)};
\node[rotate=270, xshift = -5] at (5.2,24) {\Large$\textcolor{color_dark_blue}{W_4}$};
\draw[color = color_dark_blue] (6.75,28.5) to [sV] (6.75,27.7); 
\draw[color = color_dark_blue] (7.5,27) -- (7.5,27.35) -- (6.75,27.35) -- (6.75,27.7);
\draw[color = color_dark_blue] (6.75,28.5) -- (6.75,29.0);
\draw[color = color_dark_blue] plot [smooth cycle] coordinates {(6.25,29.05) (6.75,29.0) (7.25,28.95) (7.25,29.05) (6.75,29.0) (6.25,28.95)};
\node[rotate=270] at (5.8,28.1) {\Large$\textcolor{color_dark_blue}{W_5}$};
\draw[color = color_dark_blue] (16.5,24.5)   to [sV] (15.7,24.5); 
\draw[color = color_dark_blue] (15,24.5)--(15.7,24.5);
\draw[color = color_dark_blue] plot [smooth cycle] coordinates {(17.05,25) (17,24.5) (16.95,24) (17.05,24) (17,24.5) (16.95,25)};
\draw[color = color_dark_blue] (17,24.5)--(16.5,24.5);
\node[rotate=270, xshift = -5] at (16.1,25.2) {\Large$\textcolor{color_dark_blue}{W_6}$};
\draw [color_ConvD, dashed, line width=0.5mm] plot [smooth, tension=0] coordinates { (0.5,4) (9.25,4) (9.25,7.5) (12.3, 9.5) (17.5, 9.5) (17.5, 12.15) (0.5, 12.15) (0.5,4)};
\draw [color_StochD, dashed, line width=0.5mm] plot [smooth, tension=0] coordinates {(9.45,4) (9.45,7.35) (12.4,9.3) (17.5,9.3) (17.5,4) (9.45,4)};
\draw [color_StochD, dashed, line width=0.5mm] plot [smooth, tension=0] coordinates { (0.5,12.3) (3,12.3) (14,12.3) (17.5, 18) (17.5, 20.2) (5, 20.2) (5, 16) (4, 15) (0.5, 15) (0.5, 12.3)};
\draw [color_bilevel, dashed, line width=0.5mm] plot [smooth, tension=0] coordinates { (0.5,15.2) (4,15.2) (4.8,16) (4.8,20.4) (17.5,20.4) (17.5,26) (12,29.5) (0.5,29.5) (0.5,15.2)};
\node at (10,24) {\LARGE \textcolor{color_bilevel}{Control zone III}};
\node at (10,19.5) {\LARGE \textcolor{color_StochD}{Control zone II}};
\node[rotate=270] at (1.5,9.5) {\LARGE \textcolor{color_ConvD}{Control zone I}};
\end{circuitikz}
\end{turn}
}
\caption{IEEE 24-Bus reliability test system layout with three reserve control zones.}
\label{Splitting}
\squeezeup
\end{figure}

We now consider the optimal reserve dimensioning in a multi-zone setting. For this purpose, the IEEE 24-Bus system is split into three reserve control zones as depicted in Fig. \ref{Splitting}. This zonal layout corresponds to the one proposed in \cite{Jensen_2017}. In each control zone there are at least one wind power unit with capacity of 100 MW and at least two flexible generation units. Unlike in the previous instance, the requirements computed with the probabilistic approach are now set for each reserve control zone independently considering the distribution of wind power production of each zone. The reliability level $\xi$ is set to $0.98$.

\begin{figure}[]
\centering
\resizebox{8.75cm}{!}{%
\begin{tikzpicture}[thick,scale=1]
\pgfplotsset{
every non boxed y axis/.append style={y axis line style=-},
every non boxed x axis/.append style={x axis line style=-},
every tick label/.append style={font=\tiny},
compat=newest
}

\begin{axis}[
legend pos = north west,
legend style={fill=none, draw=none, font=\scriptsize, legend cell align={left}},
ylabel near ticks,
label style={font=\scriptsize},
ylabel style={align=center},
ylabel= {$D_{z}^{\text{U}}$ [MW]},
title={\footnotesize Zone I},
solid,
xmajorticks=false,
width=4cm,
height=3.5cm,
ymax = 200,
ymin  = 0,
xmin=1,
xmax=24,
]
\addplot [color_ConvD!120, line width=0.2mm,densely dashed,smooth] table [x index = 0, y index = 1] {plot_data/reserve_requirements_three_zones.dat};
\addlegendentry{Probabilistic};
\addplot [color_bilevel!120, line width=0.2mm,solid,smooth] table [x index = 0, y index = 7] {plot_data/reserve_requirements_three_zones.dat};
\addlegendentry{Enhanced};
\end{axis}

\begin{axis}[
title={\footnotesize Zone II},
solid,
xmajorticks=false,
yticklabels={},
width=4cm,
height=3.5cm,
xshift = 75,
ymax = 200,
ymin  = 0,
xmin=1,
xmax=24,
]
\addplot [color_ConvD!120, line width=0.2mm,densely dashed,smooth] table [x index = 0, y index = 2] {plot_data/reserve_requirements_three_zones.dat};
\addplot [color_bilevel!120, line width=0.2mm,solid,smooth] table [x index = 0, y index = 8] {plot_data/reserve_requirements_three_zones.dat};
\end{axis}

\begin{axis}[
title={\footnotesize Zone III},
solid,
xmajorticks=false,
yticklabels={},
width=4cm,
height=3.5cm,
xshift = 150,
ymax = 200,
ymin  = 0,
xmin=1,
xmax=24,
]
\addplot [color_ConvD!120, line width=0.2mm,densely dashed,smooth] table [x index = 0, y index = 3] {plot_data/reserve_requirements_three_zones.dat};
\addplot [color_bilevel!120, line width=0.2mm,solid,smooth] table [x index = 0, y index = 9] {plot_data/reserve_requirements_three_zones.dat};
\end{axis}

\begin{axis}[
ylabel near ticks,
label style={font=\scriptsize},
ylabel style={align=center},
ylabel= {$D_{z}^{\text{D}}$ [MW]},
solid,
width=4cm,
height=3.5cm,
yshift = -75,
ymax = 90,
ymin  = -0.2,
xmin=1,
xmax=24,
xtick={1,6,12,18,24}
]
\addplot [color_ConvD!120, line width=0.2mm,densely dashed,smooth] table [x index = 0, y index = 4] {plot_data/reserve_requirements_three_zones.dat};
\addplot [color_bilevel!120, line width=0.2mm,solid,smooth] table [x index = 0, y index = 10] {plot_data/reserve_requirements_three_zones.dat};
\end{axis}

\begin{axis}[
xlabel={\scriptsize Operating hour},
xlabel near ticks,
solid,
yticklabels={},
width=4cm,
height=3.5cm,
xshift = 75,
yshift = -75,
ymax = 90,
ymin  = -0.2,
xmin=1,
xmax=24,
xtick={1,6,12,18,24}
]
\addplot [color_ConvD!120, line width=0.2mm,densely dashed,smooth] table [x index = 0, y index = 5] {plot_data/reserve_requirements_three_zones.dat};
\addplot [color_bilevel!120, line width=0.2mm,solid,smooth] table [x index = 0, y index = 11] {plot_data/reserve_requirements_three_zones.dat};
\end{axis}

\begin{axis}[
solid,
yticklabels={},
width=4cm,
height=3.5cm,
xshift = 150,
yshift = -75,
ymax = 90,
ymin  = -0.2,
xmin=1,
xmax=24,
xtick={1,6,12,18,24}
]
\addplot [color_ConvD!120, line width=0.2mm,densely dashed,smooth] table [x index = 0, y index = 6] {plot_data/reserve_requirements_three_zones.dat};
\addplot [color_bilevel!120, line width=0.2mm,solid,smooth] table [x index = 0, y index = 12] {plot_data/reserve_requirements_three_zones.dat};
\end{axis}
\end{tikzpicture}
}
\caption{24-hour profiles of probabilistic and enhanced reserve requirements in three control zones.}
\label{Reserve_allocation_zones}
\end{figure}
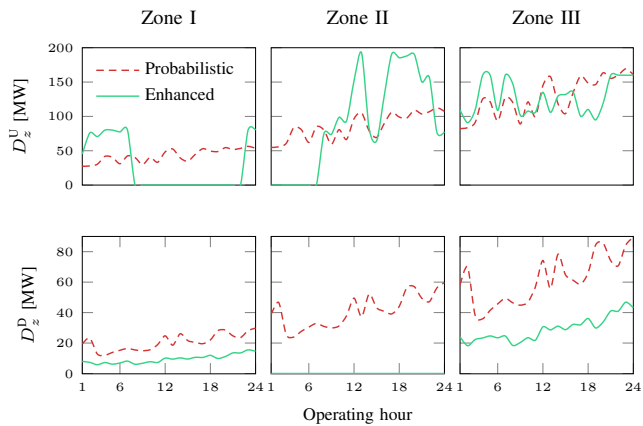

The resulting allocation for upward and downward reserve requirements among control zones is summarized in Fig. \ref{Reserve_allocation_zones}, indicating that the probabilistic approach sets the reserve needs proportionally to the amount of stochastic in-feed in the respective control zone. On the other hand, the proposed model defines the requirements considering not only the zonal wind power in-feed, but also the cost implications of procuring reserve in a specific zone. As a result, the model finds it more efficient to constantly procure upward reserve from the third zone and obtain the remaining upward reserve that is needed either from the first or the second zone depending on the operating hour. In addition, this reserve allocation indicates that it is never optimal to procure downward reserve from the second zone in terms of expected system cost.

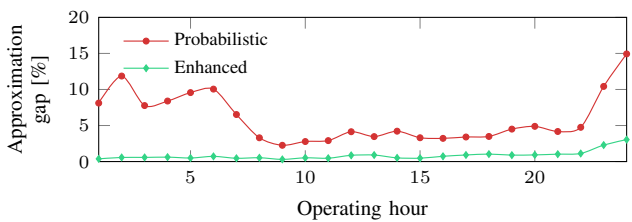
\begin{figure}[]
\center
\begin{tikzpicture}[thick,scale=1]
\pgfplotsset{ymin=0, ymax=20, xmin=1, xmax=24, try min ticks=5}
\begin{axis}[
xlabel near ticks,
ylabel near ticks,
xlabel={Operating hour},
label style={font=\footnotesize},
ylabel style={align=center},
ylabel= {Approximation \\ gap [\%]},
tick label style={font=\scriptsize},
legend pos = north west, xshift = 10,
legend style={fill=none, draw=none, name = ConvD_leg_pos, font=\scriptsize, legend cell align={left}}, width=8.6cm,
height=3.5cm
]
\addplot [color_ConvD!120, line width=0.15mm, mark=*, mark size=1.1, smooth, name path=A] table [x index = 0, y index = 1] {plot_data/efficiency_gap.dat};
\addlegendentry{Probabilistic}
\addplot [color_bilevel!120, line width=0.15mm, mark=diamond*, mark size=1.25, smooth] table [x index = 0, y index = 2] {plot_data/efficiency_gap.dat};
\addlegendentry{Enhanced}
\end{axis}
\end{tikzpicture}
\caption{Approximation gap of the sequential market with probabilistic and enhanced reserve requirements compared to stochastic dispatch.}
\label{three_zone_cost}
\squeezeup
\end{figure}

This optimal reserve allocation among control zones is supported by the approximation gap depicted in Fig. \ref{three_zone_cost}, \textcolor{WildStrawberry}{showing the relative cost difference of the sequential market with respect to the ideal solution.} The requirements provided by the proposed model efficiently approximate the ideal solution with nearly zero gap over the first operating hours, and this gap remains relatively small for the subsequent hours \textcolor{WildStrawberry}{as opposed to the large gap when probabilistic requirements are used}. \textcolor{ForestGreen}{The definition of multiple control zones allows to set enhanced reserve requirements that are closer to the `generator-specific' reserve allocation of the stochastic model. } Indeed, \textcolor{ForestGreen}{compared to the single-zone setup} in section \ref{Impact_of_reserve_requirements}, the operating cost \textcolor{ForestGreen}{reduces by 2.5\%, from \$24,408 to \$24,034, after the definition of three control zones}.

\vspace{-0.03cm}
\subsection{Assessing enhanced reserve requirements in the presence of non-convexities} \label{24RTS_with_UC_conctrains}

\textcolor{blue}{To assess the performance for the proposed reserve quantification model as a proxy model for the power markets with a more comprehensive and non-convex representation of technical constraints, we use the enhanced reserve requirements provided by the proposed model \eqref{prob:bilevel_clearing} as inputs to the sequential market-clearing problem \eqref{prob:reserve_clearing}-\eqref{prob:real_time_clearing} with unit commitment and ramping constraints integrated in the day-ahead auction as explained in Appendix \ref{appA}.}

\textcolor{blue}{Figure \ref{Cost_proxy} shows the hourly profile of expected operating system cost resulting from the implementation of the enhanced requirements in the system with full representation of the technical constraints. This profile is compared against those obtained by setting probabilistic reserve requirements with reliability levels of 98\% and 90\%. The reserve requirements provided by the proposed model always attain better cost efficiency than the probabilistic requirements, even though the proposed model does not account for the whole set of technical limits of power plants. In the first case in Fig.\ref{Cost_proxy} (a), the model allows savings of \$23,746 that nearly equal to the cost of peak-hour operation, and it allows even larger savings of \$28,845 in the second case in Fig.\ref{Cost_proxy} (b).}

\begin{center}
\begin{figure}[]
\centering
\subfigure[Reliability level $\xi = 98\%$]
{
\resizebox{8cm}{!}{%
\begin{tikzpicture}[thick,scale=1]
\begin{axis}[
ymin=4,
ymax=28,
xmax=24,
xmin=1,
xtick={1,3,...,24},
ytick={4,10,...,26},
try min ticks=10,
max space between ticks=10pt,
xlabel={Operating hour},
ylabel={Expected cost [$\$10^3$]} ,
x label style = {at={(axis description cs:0.5,-0.18)},anchor=south, font=\scriptsize},
y label style=  {at={(axis description cs:0.10,.5)},rotate=0,anchor=south, font=\scriptsize},
legend style={draw=none, fill=none, legend columns=1, font=\scriptsize, legend pos={north west}},
tick label style={font=\scriptsize},
width=8.5cm,
height=4cm,
]
\addplot [draw=color_ConvD, name path=prob] table [x index = 0, y index = 1] {plot_data/proxy_pr_1per.dat};
\addlegendentry{Probabilistic solution};
\addplot [draw=color_bilevel, name path=enh] table [x index = 0, y index = 1] {plot_data/proxy_enahnced.dat};
\addlegendentry{Enhanced solution};
\addplot [thick,color=blue,fill=blue,fill opacity=0.05] fill between[of=prob and enh];
\addplot[mark=o, fill = blue!50, mark size = 0.8pt, color=blue!50, line width = 0.15mm] coordinates {(19.5,21)} node[pin={[pin distance=0.4cm]-93:{\textcolor{black}{\scriptsize Cumulative saving of \$23,746}}}]{};
\end{axis}
\end{tikzpicture}
}
}
\\
\subfigure[Reliability level $\xi = 90\%$]
{
\resizebox{8cm}{!}{%
\begin{tikzpicture}[thick,scale=1]
\begin{axis}[
ymin=4,
ymax=28,
xmax=24,
xmin=1,
xtick={1,3,...,24},
ytick={4,10,...,26},
try min ticks=10,
max space between ticks=10pt,
xlabel={Operating hour},
ylabel={Expected cost [$\$10^3$]} ,
x label style = {at={(axis description cs:0.5,-0.18)},anchor=south, font=\scriptsize},
y label style=  {at={(axis description cs:0.10,.5)},rotate=0,anchor=south, font=\scriptsize},
legend style={draw=none, fill=none, legend columns=-1},
tick label style={font=\scriptsize},
width=8.5cm,
height=4cm,
]
\addplot [draw=color_ConvD, name path=prob] table [x index = 0, y index = 1] {plot_data/proxy_pr_5per.dat};
\addplot [draw=color_bilevel, name path=enh] table [x index = 0, y index = 1] {plot_data/proxy_enahnced.dat};
\addplot [thick,color=blue,fill=blue,fill opacity=0.05] fill between[of=prob and enh];
\addplot[mark=o, fill = blue!50, mark size = 0.8pt, color=blue!50, line width = 0.15mm] coordinates {(19.5,21.8)} node[pin={[pin distance=0.6cm]-93:{\textcolor{black}{\scriptsize Cumulative saving of \$28,845}}}]{};
\end{axis}
\end{tikzpicture}
}
}
\caption{Expected operating cost yielded by the implementation of the probabilistic and enhanced reserve requirements in the conventional market-clearing problem  \eqref{prob:reserve_clearing}-\eqref{prob:real_time_clearing} including the unit commitment constraints \eqref{app1}-\eqref{app7}.}
\label{Cost_proxy}
\squeezeup
\end{figure}
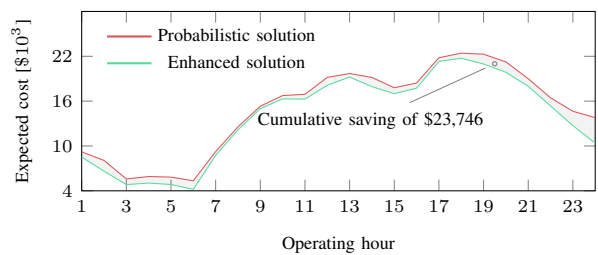
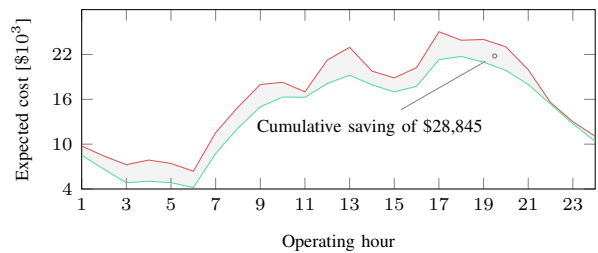
\end{center}

\vspace{-1.3cm}
\textcolor{blue}{
\subsection{Application to the IEEE-96 RTS} \label{IEEE96study}
We now consider the modernized version of the IEEE-96 RTS Test System proposed in \cite{IEEE96RTS} to assess the scalability of the proposed model. The test system includes three control zones interconnected by six tie-lines. The system demand follows a 24-hour profile with a peak load of 7.5 GW. The conventional generation is represented by 6 nuclear power plants serving the base load, 3 coal power plants that offer 40\% of their capacities for the reserve needs, and 87 gas-fired power plants offering 100\% of their capacities to the reserve procurement auction. The reserve offering prices of flexible units are set to 25\% of marginal production cost for both upward and downward reserve needs. There are 19 wind farms distributed among the control zones with the overall capacity of 2.76 GW. Their stochastic output is described by  100 equiprobable scenarios obtained from \cite{pinson2013wind}. The input data and the corresponding GAMS codes are provided in the electronic companion of the paper \cite{companion}.
}

\textcolor{blue}{
The test case is solved for wind penetration levels of 13.8\%, 23.0\%, and 36.8\% of the peak-hour load by implementing the multicut Bender's algorithm explained in Section \ref{section_sol_str}. The tolerance of the algorithm is set to 0.02\% requiring three to eight iterations depending on the operating hour. The resulting CPU time is reported in Table \ref{CPU_time}. The CPU time in all three cases is kept below one hour allowing timely day-ahead planning with the proposed model. It is worth mentioning that the CPU time can be reduced at the expense of a marginal deviation from the global optimum with higher tolerance.
}
\textcolor{blue}{
\begin{table}[]
\centering
\textcolor{blue}{
\caption{\textcolor{blue}{CPU performance of the Bender's algorithm.}}
\renewcommand{\arraystretch}{1.5}
\begin{tabular}{c|ccc}
\hline
Wind penetration {[}\%{]} & 13.8 & 23.0 & 36.8 \\
CPU time {[}min{]} & 32.1 & 33.5 & 58.6 \\
\hline
\end{tabular}
\label{CPU_time}
}
\end{table}
}

\textcolor{blue}{
The daily operating cost resulting from the implementation of the enhanced zonal reserve requirements computed by the proposed model is always lower than those provided by the probabilistic approach with reliability levels of 90\% and 98\%, as demonstrated in Table \ref{IEEE96_cost}. The difference in operating cost is explained by the anticipated cost of procuring upward and downward reserves from a specific control zone, while the probabilistic requirements are solely obtained proportionally to the amount of stochastic in-feed in control zones. As a result, the relative cost savings provided by the model increases with the wind penetration level and ranges between 0.6\% and 7.2\%. Further cost savings towards the ideal solution provided by the stochastic model is limited due to the enforced merit order in both reserve and day-ahead markets.}
\begin{table}[]
\centering
\textcolor{blue}{
\caption{\textcolor{blue}{Daily operating cost with probabilistic and enhanced zonal reserve requirements in comparison with the stochastic ideal solution [\$1000].}}
\renewcommand{\arraystretch}{1}
\begin{tabular}{ccccc}
\hline
\multirow{2}{*}{\begin{tabular}[c]{@{}c@{}}Wind penetration\\ {[}\%{]}\end{tabular}} & \multicolumn{2}{c}{\begin{tabular}[c]{@{}c@{}}Probabilistic \\ solution\end{tabular}} & \multirow{2}{*}{\begin{tabular}[c]{@{}c@{}}Enhanced\\ solution\end{tabular}} &
\multirow{2}{*}{\begin{tabular}[c]{@{}c@{}}Ideal\\ solution\end{tabular}}
\\
\cline{2-3}
 & $\xi=90\%$ & $\xi=98\%$ & & \\
\hline
13.8 & 1,912.4 & 1,888.8 & 1,877.3 & 1,850.0\\
23.0 & 1,760.8 & 1,719,3 &  1,700.8 & 1,660.5\\
36.8 & 1,550.7 & 1,482.3 & 1,446.0 & 1,402.8\\
\hline
\end{tabular}
\label{IEEE96_cost}
}
\squeezeup
\end{table}
\textcolor{blue}{
Finally, Table \ref{IEEE96_cost2} illustrates the economic benefit that the proposed model yields as a proxy for the system with the full network representation and technical constraints of power plants described in Appendix \ref{appA}. The results show that in spite of the incomplete description of technical constraints in the lower level of the proposed bilevel model, it still provides a feasible input with a sensible cost reduction for the markets with non-convexities. The economic benefit provided by the model ranges from 0.5\% to 1.6\%. Moreover, the proposed approach further outperforms the probabilistic one for the largest wind penetration level, where the overestimated requirements provided by the probabilistic approach lead to a reserve schedule that results in an infeasible day-ahead operation.
}

\begin{table}[]
\centering
\textcolor{blue}{
\caption{\textcolor{blue}{Daily operating cost with probabilistic and enhanced zonal reserve requirements with full representation of technical constraints [\$1000].}}
\renewcommand{\arraystretch}{1}
\begin{tabular}{cccc}
\hline
\multirow{2}{*}{\begin{tabular}[c]{@{}c@{}}Wind penetration\\ {[}\%{]}\end{tabular}} & \multicolumn{2}{c}{\begin{tabular}[c]{@{}c@{}}Probabilistic \\ solution\end{tabular}} & \multirow{2}{*}{\begin{tabular}[c]{@{}c@{}}Enhanced\\ solution\end{tabular}} \\
\cline{2-3}
 & $\xi=90\%$ & $\xi=98\%$ &  \\
\hline
13.8 & 2,072.2 & 2,073.4 &  2,061.5 \\
23.0 & 1,947.9 & 1,949.1 & 1,928.6  \\
36.8 & 1,764.4 & infeas. & 1,735.9 \\
\hline
\end{tabular}
\squeezeup
\label{IEEE96_cost2}
}
\end{table}

\section{Conclusion} \label{sec4}

This paper considers the optimal setting of reserve requirements in a European market framework. We propose a new method to quantify reserve needs that brings \textcolor{blue}{the sequence of the reserve, day-ahead and real-time markets} closer to the ideal stochastic energy and reserves co-optimization model in terms of total expected cost. The proposed model is formulated as a stochastic bilevel problem, \textcolor{blue}{which is eventually recast as a MILP problem. To reduce the computational burden of this model, we apply an iterative solution approach based on the multi-cut Bender's decomposition algorithm.}

\textcolor{blue}{Our numerical studies demonstrate the benefit of properly setting reserve requirements.
Our reserve quantification model outperforms both the probabilistic and the stochastic reserve setting approaches due to its preemptive ability to anticipate the impact of day-ahead decisions on the real-time operation, while taking into account the actual market structure.}
\textcolor{blue}{Considering the increasing penetration of stochastic power producers, we show that the reserve requirements provided by the proposed model take the expected system operating cost closer to that given by the ideal energy and reserve co-optimization model, but the degree of this approximation is limited due to the sequential scheduling of reserve and energy in European electricity markets.}
\textcolor{blue}{However, our analysis further indicates that the definition of multiple reserve control zones allows for a more efficient spatial allocation of reserves, which reduces the approximation gap with respect to the ideal stochastic model.}
\textcolor{blue}{Finally, the efficiency of the proposed reserve dimensioning model was tested against market designs whose clearing process explicitly account for inter-temporal and non-convex constraints, i.e. ramping limits and unit commitment constraints. Even though the proposed model does not account for the whole set of technical constraints of such markets, the enhanced reserve requirements still bring the cost of sequential market operation closer to the stochastic ideal, highlighting the importance of the intertemporal coordination between the three trading floors through the intelligent setting of reserve needs.}

\textcolor{blue}{Future research may focus on the consideration of the tight relaxations of the unit commitment constraints to achieve better approximations for the case of non-convex market designs, and the corresponding tuning of the Bender's decomposition algorithm to better cope with the intertemporal constraints.}

\appendix

\subsection{Incorporation of unit commitment and ramping constraints} \label{appA}

\textcolor{blue}{In contrast to the prevailing approach of the European market design, other electricity markets, e.g., the majority of US markets, explicitly model unit commitment constraints and thermal limits of power plants in the market-clearing problem. To assess the performance of the proposed reserve quantification model in markets with unit commitment constraints, the following set of constraints are integrated in the day-ahead market-clearing problem:}
\textcolor{blue}{
\begingroup
\allowdisplaybreaks
\begin{subequations}
\begin{align}
& u_{it} \underline{P}_{i} \leq P_{it}^{\text{C}} \leq u_{it} \overline{P}_{i},
\; \forall i \in I,
\; \forall t \in T,
\label{app1}\\
& SU_{it} \geq C_{i}^{\text{SU}} (u_{it} - u_{i(t-1)}),
\; \forall i \in I,
\; \forall t>1,
\label{app2}\\
& SU_{it} \geq C_{i}^{\text{SU}} (u_{it} - u_{i}^{0}),
\; \forall i \in I,
\;  t=1,
\label{app3}\\
& P_{it}^{\text{C}} - P_{i(t-1)}^{\text{C}} \leq R_{i}^{+},
\; \forall i \in I,
\; \forall t >1,
\label{app4}\\
& P_{it}^{\text{C}} - P_{i}^{\text{C},0} \leq R_{i}^{+},
\; \forall i \in I,
\;  t = 1,
\label{app5}\\
& P_{i(t-1)}^{\text{C}} - P_{it}^{\text{C}} \leq R_{i}^{-},
\; \forall i \in I,
\; \forall t >1,
\label{app6}\\
& P_{i}^{\text{C},0} - P_{it}^{\text{C}} \leq R_{i}^{-},
\; \forall i \in I,
\;  t =1,
\label{app7}
\end{align}
\end{subequations}
\endgroup
where $t\in T$ is the set of operating hours, $C_{i}^{\text{SU}}$ is a start-up cost of unit $i$, $R_{i}^{+}$ and $R_{i}^{-}$ are the ramp-up and ramp-down limits, $\underline{P}_{i}$ is a minimum power output limit, and $P_{i}^{\text{C},0}$ and $u_{i}^{0}$ are the initial power output and commitment status of unit $i$. The set of decision variables of the original problem is supplemented with variable $u_{it}\in \{0,1\}$ that denotes the commitment status of generating units, and variable $SU_{it}$ that computes the cost induced by the start-up of generating units. Now, the generating limits of each unit are additionally enforced by commitment decisions of the system operator by \eqref{app1}. Binary logic is controlled by \eqref{app2} and \eqref{app3} and activated by augmenting $SU_{it}$ into the original objective function of problem \eqref{prob:day_ahead_clearing}. The ramp limits of generators are accounted for through \eqref{app4}-\eqref{app7}.}

\bibliographystyle{IEEEtran}
\bibliography{references}

\begin{thebibliography}{10}
\providecommand{\url}[1]{#1}
\csname url@samestyle\endcsname
\providecommand{\newblock}{\relax}
\providecommand{\bibinfo}[2]{#2}
\providecommand{\BIBentrySTDinterwordspacing}{\spaceskip=0pt\relax}
\providecommand{\BIBentryALTinterwordstretchfactor}{4}
\providecommand{\BIBentryALTinterwordspacing}{\spaceskip=\fontdimen2\font plus
\BIBentryALTinterwordstretchfactor\fontdimen3\font minus
  \fontdimen4\font\relax}
\providecommand{\BIBforeignlanguage}[2]{{%
\expandafter\ifx\csname l@#1\endcsname\relax
\typeout{** WARNING: IEEEtran.bst: No hyphenation pattern has been}%
\typeout{** loaded for the language `#1'. Using the pattern for}%
\typeout{** the default language instead.}%
\else
\language=\csname l@#1\endcsname
\fi
#2}}
\providecommand{\BIBdecl}{\relax}
\BIBdecl

\bibitem{Aigner_2012}
T.~Aigner, S.~Jaehnert, G.~L. Doorman, and T.~Gjengedal, ``The effect of
  large-scale wind power on system balancing in {Northern Europe},'' \emph{IEEE
  Trans. Sustain. Energy}, vol.~3, no.~4, pp. 751--759, 2012.

\bibitem{Morales_2012}
J.~M. Morales, A.~J. Conejo, K.~Liu, and J.~Zhong, ``Pricing electricity in
  pools with wind producers,'' \emph{IEEE Trans. Power Syst.}, vol.~27, no.~3,
  pp. 1366--1376, 2012.

\bibitem{1525111}
F.~Bouffard, F.~D. Galiana, and A.~J. Conejo, ``Market-clearing with stochastic
  security-part {I}: formulation,'' \emph{IEEE Trans. Power Syst.}, vol.~20,
  no.~4, pp. 1818--1826, 2005.

\bibitem{Papavasiliou_2015}
A.~Papavasiliou, S.~S. Oren, and B.~Rountree, ``Applying high performance
  computing to transmission-constrained stochastic unit commitment for
  renewable energy integration,'' \emph{IEEE Trans. Power Syst.}, vol.~30,
  no.~3, pp. 1109--1120, 2015.

\bibitem{4806110}
A.~Tuohy, P.~Meibom, E.~Denny, and M.~O'Malley, ``Unit commitment for systems
  with significant wind penetration,'' \emph{IEEE Trans. Power Syst.}, vol.~24,
  no.~2, pp. 592--601, 2009.

\bibitem{4556639}
J.~Wang, M.~Shahidehpour, and Z.~Li, ``Security-constrained unit commitment
  with volatile wind power generation,'' \emph{IEEE Trans. Power Syst.},
  vol.~23, no.~3, pp. 1319--1327, 2008.

\bibitem{zavala2017stochastic}
V.~M. Zavala, K.~Kim, M.~Anitescu, and J.~Birge, ``A stochastic electricity
  market clearing formulation with consistent pricing properties,'' \emph{Oper.
  Res.}, vol.~65, no.~3, pp. 557--576, 2017.

\bibitem{morales2014electricity}
J.~M. Morales, M.~Zugno, S.~Pineda, and P.~Pinson, ``Electricity market
  clearing with improved scheduling of stochastic production,'' \emph{Eur. J.
  Oper. Res.}, vol. 235, no.~3, pp. 765--774, 2014.

\bibitem{Kazempour_2018}
J.~Kazempour, P.~Pinson, and B.~F. Hobbs, ``A stochastic market design with
  revenue adequacy and cost recovery by scenario: Benefits and costs,''
  \emph{IEEE Trans. Power Syst.}, vol.~33, no.~4, pp. 3531--3545, 2018.

\bibitem{SARFATI2018851}
M.~Sarfati, M.~R. Hesamzadeh, D.~R. Biggar, and R.~Baldick, ``Probabilistic
  pricing of ramp service in power systems with wind and solar generation,''
  \emph{Renewable and Sustainable Energy Reviews}, vol.~90, pp. 851 -- 862,
  2018.

\bibitem{Jensen_2017}
T.~V. Jensen, J.~Kazempour, and P.~Pinson, ``Cost-optimal {ATCs} in zonal
  electricity markets,'' \emph{IEEE Trans. Power Syst.}, vol.~33, no.~4, pp.
  3624--3633, 2018.

\bibitem{Wang_2013}
B.~Wang and B.~F. Hobbs, ``Flexiramp market design for real-time operations:
  Can it approach the stochastic optimization ideal?'' in \emph{2013 IEEE
  PESGM}, 2013, pp. 1--5.

\bibitem{hogan2013electricity}
W.~W. Hogan, ``Electricity scarcity pricing through operating reserves,''
  \emph{EEEP}, vol.~2, no.~2, pp. 65--86, 2013.

\bibitem{papavasiliou2017remuneration}
A.~Papavasiliou and Y.~Smeers, ``Remuneration of flexibility using operating
  reserve demand curves: A case study of {Belgium}.'' \emph{Energy Journal},
  vol.~38, no.~6, 2017.

\bibitem{rebours2005survey}
Y.~Rebours and D.~Kirschen, ``A survey of definitions and specifications of
  reserve services,'' Tech. Rep., 2005.

\bibitem{manual2012energy}
\BIBentryALTinterwordspacing
``{Emergency Operations},'' {\textit{PJM Manual}}, Tech. Rep., 2011. [Online].
  Available:
  \url{http://www.pjm.com/~/media/documents/manuals/archive/m13/m13v46-emergency-operations-11-16-2011.ashx}
\BIBentrySTDinterwordspacing

\bibitem{strbac2007impact}
G.~Strbac, A.~Shakoor, M.~Black, D.~Pudjianto, and T.~Bopp, ``Impact of wind
  generation on the operation and development of the {UK} electricity
  systems,'' \emph{Electr. Pow. Syst. Res.}, vol.~77, no.~9, pp. 1214--1227,
  2007.

\bibitem{lee2012analyzing}
D.~Lee and R.~Baldick, ``Analyzing the variability of wind power output through
  the power spectral density,'' in \emph{2012 IEEE PESGM}, 2012, pp. 1--8.

\bibitem{Doherty_1425549}
R.~Doherty and M.~O'Malley, ``A new approach to quantify reserve demand in
  systems with significant installed wind capacity,'' \emph{IEEE Trans. Power
  Syst.}, vol.~20, no.~2, pp. 587--595, 2005.

\bibitem{5929570}
Y.~V. Makarov, P.~V. Etingov, J.~Ma, Z.~Huang, and K.~Subbarao, ``Incorporating
  uncertainty of wind power generation forecast into power system operation,
  dispatch, and unit commitment procedures,'' \emph{IEEE Trans. Sustain.
  Energy}, vol.~2, no.~4, pp. 433--442, 2011.

\bibitem{6942382}
Y.~Dvorkin, D.~S. Kirschen, and M.~A. Ortega-Vazquez, ``Assessing flexibility
  requirements in power systems,'' \emph{{IET} GENER. TRANSM. DIS.}, vol.~8,
  no.~11, pp. 1820--1830, 2014.

\bibitem{7084167}
H.~Nosair and F.~Bouffard, ``Flexibility envelopes for power system operational
  planning,'' \emph{IEEE Trans. Sustain. Energy}, vol.~6, no.~3, pp. 800--809,
  2015.

\bibitem{7552596}
------, ``Economic dispatch under uncertainty: The probabilistic envelopes
  approach,'' \emph{IEEE Trans. Power Syst.}, vol.~32, no.~3, pp. 1701--1710,
  2017.

\bibitem{lange2005uncertainty}
M.~Lange, ``On the uncertainty of wind power predictions—analysis of the
  forecast accuracy and statistical distribution of errors,'' \emph{J. Sol.
  Energy Eng.}, vol. 127, no.~2, pp. 177--184, 2005.

\bibitem{5565529}
M.~A. Matos and R.~J. Bessa, ``Setting the operating reserve using
  probabilistic wind power forecasts,'' \emph{IEEE Trans. Power Syst.},
  vol.~26, no.~2, pp. 594--603, 2011.

\bibitem{6299425}
H.~Holttinen, M.~Milligan, E.~Ela, N.~Menemenlis, J.~Dobschinski, B.~Rawn,
  R.~J. Bessa, D.~Flynn, E.~Gomez-Lazaro, and N.~K. Detlefsen, ``Methodologies
  to determine operating reserves due to increased wind power,'' \emph{IEEE
  Trans. Sustain. Energy}, vol.~3, no.~4, pp. 713--723, 2012.

\bibitem{6487424}
P.~N. Biskas, D.~I. Chatzigiannis, and A.~G. Bakirtzis, ``European electricity
  market integration with mixed market designs--part {I}: Formulation,''
  \emph{IEEE Trans. Power Syst.}, vol.~29, no.~1, pp. 458--465, 2014.

\bibitem{Kasina_2014}
S.~Kasina, ``Essays on unit commitment and interregional cooperation in
  transmission planning,'' 2017, {Ph.D.} thesis.

\bibitem{7914790}
J.~Kazempour and B.~F. Hobbs, ``Value of flexible resources, virtual bidding,
  and self-scheduling in two-settlement electricity markets with wind
  generation—-part {I}: Principles and competitive model,'' \emph{IEEE Trans.
  Power Syst.}, vol.~33, no.~1, pp. 749--759, 2018.

\bibitem{ruiz2013revealing}
C.~Ruiz, A.~J. Conejo, and D.~J. Bertsimas, ``Revealing rival marginal offer
  prices via inverse optimization,'' \emph{IEEE Trans. Power Syst.}, vol.~28,
  no.~3, pp. 3056--3064, 2013.

\bibitem{mitridati2017bayesian}
L.~Mitridati and P.~Pinson, ``A bayesian inference approach to unveil supply
  curves in electricity markets,'' \emph{IEEE Trans. Power Syst.}, vol.~33,
  no.~3, 2017.

\bibitem{pozo2017basic}
D.~Pozo, E.~Sauma, and J.~Contreras, ``Basic theoretical foundations and
  insights on bilevel models and their applications to power systems,''
  \emph{Ann. Oper. Res.}, pp. 1--32, 2017.

\bibitem{conejo2006decomposition}
A.~J. Conejo, E.~Castillo, R.~Minguez, and R.~Garcia-Bertrand,
  \emph{Decomposition techniques in mathematical programming: engineering and
  science applications}.\hskip 1em plus 0.5em minus 0.4em\relax Springer
  Science \& Business Media, 2006.

\bibitem{Ordoudis_2016}
C.~Ordoudis, P.~Pinson, J.~M. Morales, and M.~Zugno, ``An updated version of
  the {IEEE RTS 24-Bus} system for electricity market and power system
  operation studies,'' \emph{Technical University of Denmark}, 2016.

\bibitem{SWIDER20071297}
D.~J. Swider and C.~Weber, ``Bidding under price uncertainty in multi-unit
  pay-as-bid procurement auctions for power systems reserve,'' \emph{Eur. J.
  Oper. Res.}, vol. 181, no.~3, pp. 1297 -- 1308, 2007.

\bibitem{1525135}
M.~A. Plazas, A.~J. Conejo, and F.~J. Prieto, ``Multimarket optimal bidding for
  a power producer,'' \emph{IEEE Trans. Power Syst.}, vol.~20, no.~4, pp.
  2041--2050, 2005.

\bibitem{companion}
\BIBentryALTinterwordspacing
V.~Dvorkin, S.~Delikaraoglou, and J.~M. Morales. (2018) Online appendix of the
  paper ``{Setting} reserve requirements to approximate the efficiency of the
  stochastic dispatch''. [Online]. Available:
  \url{https://doi.org/10.5281/zenodo.1408881}
\BIBentrySTDinterwordspacing

\bibitem{IEEE96RTS}
\BIBentryALTinterwordspacing
H.~Pandzic, Y.~Dvorkin, T.~Qiu, Y.~Wang, and D.~Kirschen, ``Unit commitment
  under uncertainty - {GAMS} models,'' Library of the Renewable Energy Analysis
  Lab (REAL), University of Washington, Seattle, USA. [Online]. Available:
  \url{http://www.ee.washington.edu/research/real/gams_code.html}
\BIBentrySTDinterwordspacing

\bibitem{pinson2013wind}
P.~Pinson, ``Wind energy: Forecasting challenges for its operational
  management,'' \emph{Stat. Sci.}, pp. 564--585, 2013.

\end{thebibliography}

\begin{IEEEbiography}[{\includegraphics[width=1in,height=1.25in,clip,keepaspectratio]{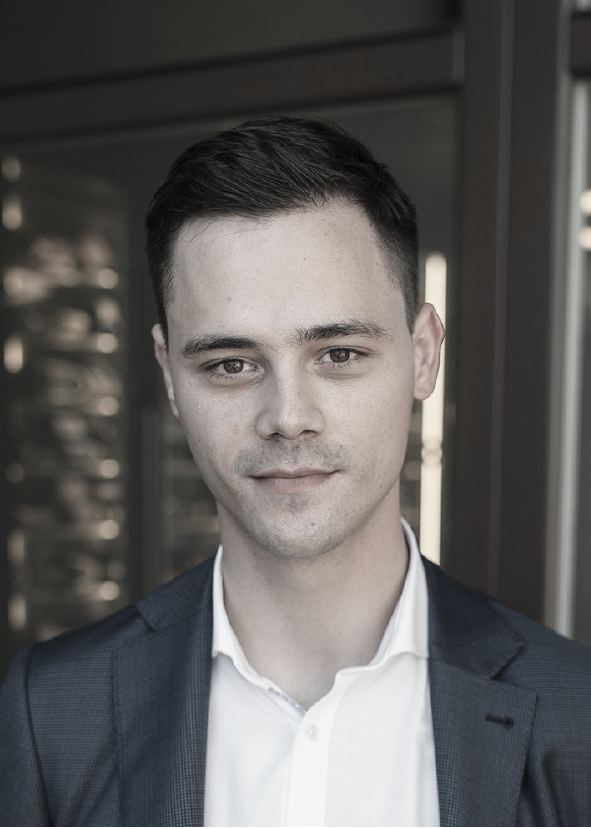}}]{Vladimir Dvorkin Jr.} (S'18) received the B.S. degree in electrical engineering from the Moscow Power Engineering Institute, Russia, in 2012, the M.Sc. degree in Economics from the Higher School of Economics, Russia, in 2014, and the M.Sc. degree in Sustainable Energy from the Technical University of Denmark in 2017. Currently he is pursuing the Ph.D. degree with the Department of Electrical Engineering,
Center for Electric Power and Energy, Technical University of Denmark.

His research interests include economics, game theory, optimization, and their applications to power systems and electricity markets.
\end{IEEEbiography}

\begin{IEEEbiography}[{\includegraphics[width=1in,height=1.25in,clip,keepaspectratio]{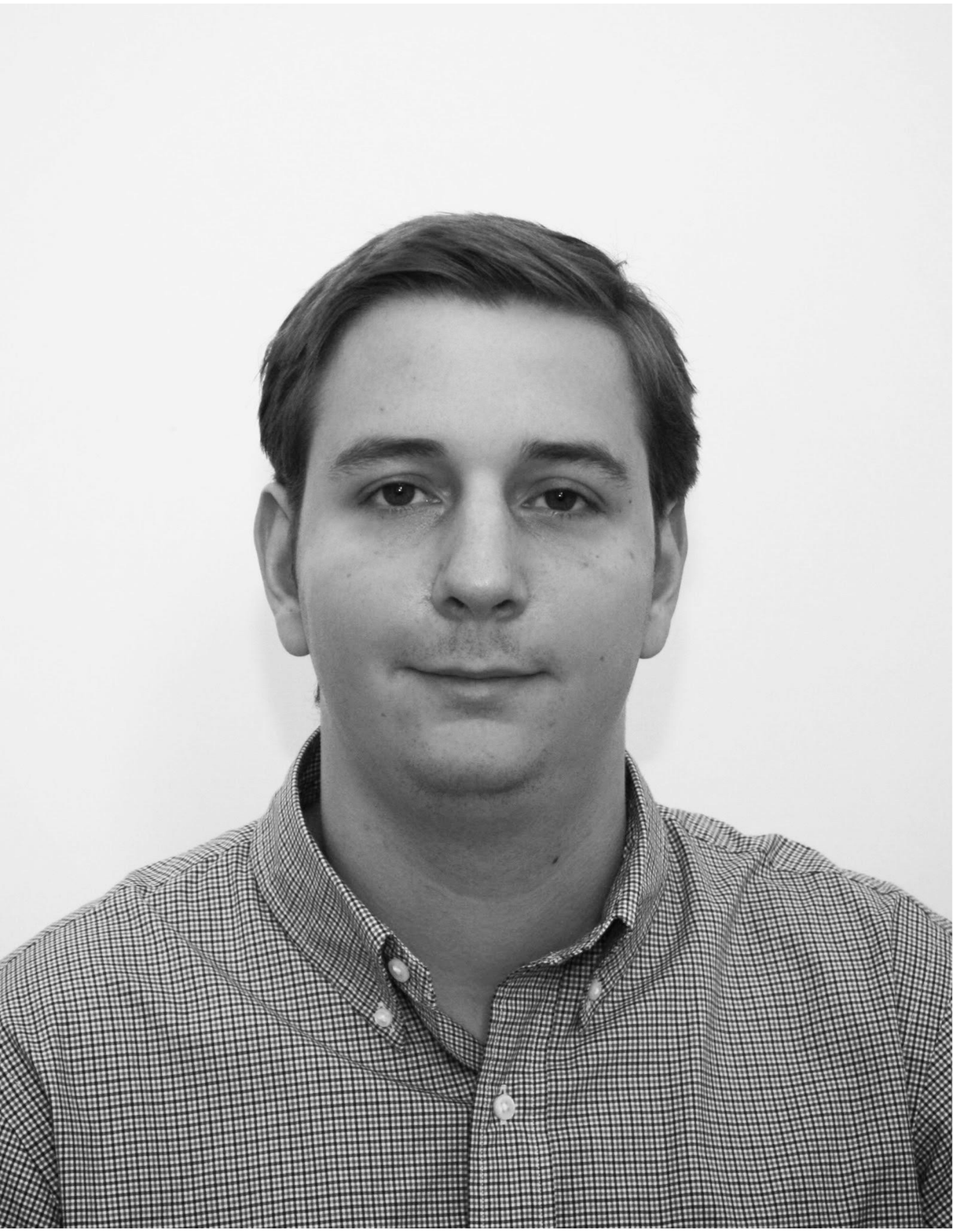}}]{Stefanos Delikaraoglou} (S'14 - M’18) received the Dipl.-Eng. degree from the School of Mechanical Engineering, National Technical University of Athens, Greece, in 2010 and the M.Sc. degree in Sustainable Energy from the Technical University of Denmark in 2012. He holds a Ph.D. degree awarded in 2016 by the Department Electrical Engineering at the Technical University of Denmark. He is currently a Postdoctoral Fellow with the EEH-Power Systems Laboratory at the Swiss Federal Institute of Technology (ETH), Zurich, Switzerland.

His research interests include energy markets and multi-energy systems modeling, decision-making under uncertainty, equilibrium models and hierarchical optimization.
\end{IEEEbiography}

\begin{IEEEbiography}[{\includegraphics[width=1in,height=1.25in,clip,keepaspectratio]{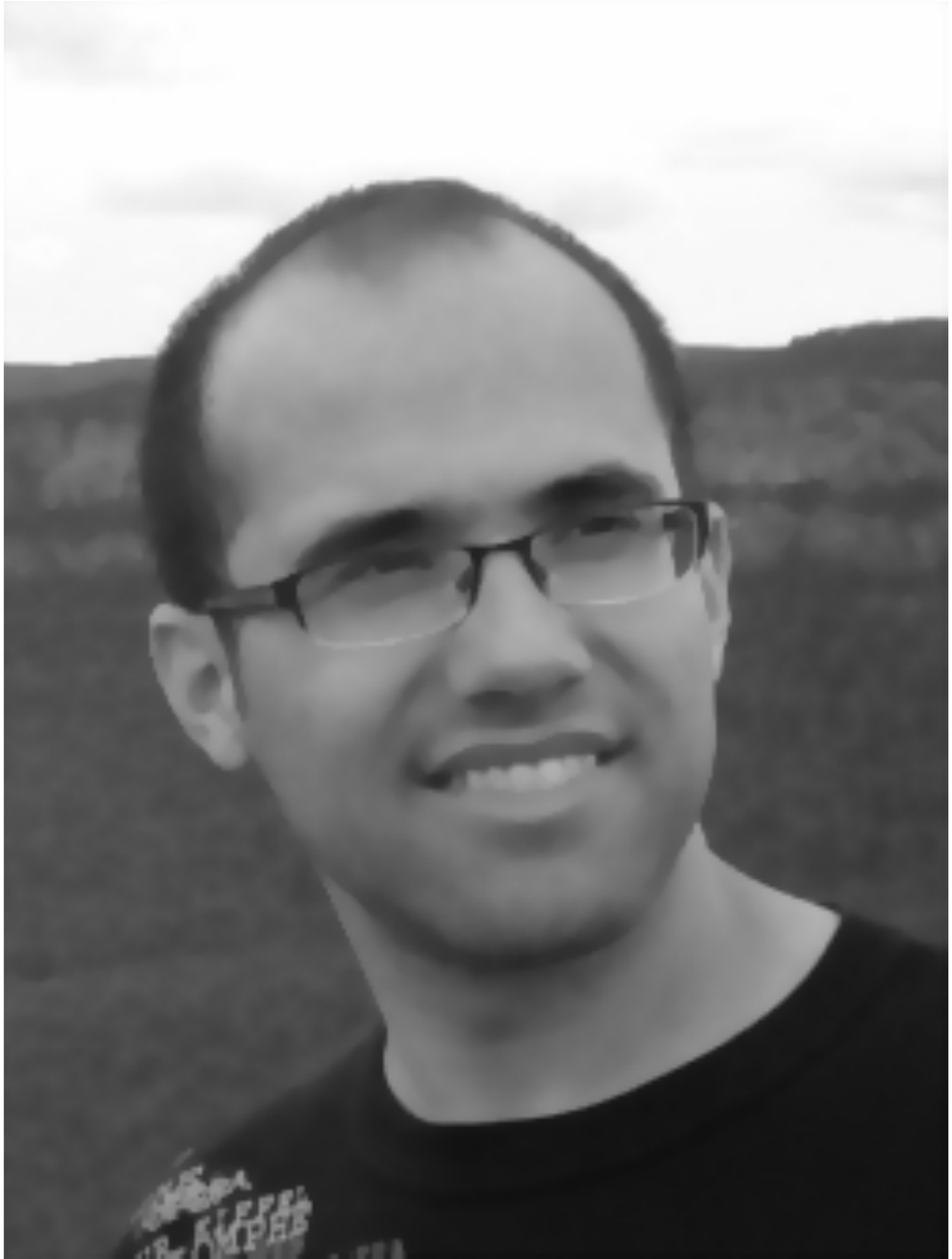}}]{Juan M.\ Morales} (S'07-M'11-SM'16) received the Ingeniero Industrial
degree from the University of  M\'alaga, M\'alaga, Spain, in 2006,
and a Ph.D. degree in Electrical Engineering from the University of Castilla-La Mancha, Ciudad Real, Spain, in 2010. He is
currently an associate professor in the Department of Applied
Mathematics at the University of M\'alaga in Spain.
\newline
\indent His research interests are in the fields of power systems
economics, operations and planning; energy analytics and optimization; smart
grids; decision-making under uncertainty, and electricity
markets.
\end{IEEEbiography}

\end{document}